\def\nosiam{}
\documentclass{article}

\input{PFKrylov.include}

\newcommand{\TheTitle}{Pipelined, Flexible Krylov Subspace Methods} 
\newcommand{\TheAuthors}{P. Sanan, S. M. Schnepp, and D. A. May}

\ifx\undefined\nosiam\headers{\TheTitle}{\TheAuthors}\fi

\title{{\TheTitle}}

\author{
P.~Sanan\footnotemark[2]\ \footnotemark[5]\ \hspace{2pt}\footnotemark[6]
\and S.M.~Schnepp\footnotemark[3]\ \footnotemark[6] 
\and D.A~May\footnotemark[4]\ \footnotemark[6]
}

\begin{document}

% **************** Uncomment this to hide all notes **************************%
%\def\hidenotes
% ****************************************************************************%

\maketitle
\renewcommand{\thefootnote}{\fnsymbol{footnote}}

\footnotetext[2]{patrick.sanan@usi.ch, patrick.sanan@erdw.ethz.ch}
\footnotetext[3]{mail@saschaschnepp.net, schnepps@ethz.ch}
\footnotetext[4]{dave.may@erdw.ethz.ch}
\footnotetext[5]{Advanced Computing Laboratory, Institute of Computational
  Science, Universit\`{a} della Svizzera italiana (USI), Via Giuseppe Buffi 13, 6904 Lugano, Switzerland}
\footnotetext[6]{Institute of Geophysics, ETH Z{\"u}rich, Sonneggstrasse 5, 8092 Z{\"u}rich, Switzerland}
\renewcommand{\thefootnote}{\arabic{footnote}}

\ifx\undefined\nosiam\slugger{sisc}{201x}{xx}{x}{x--x}\fi

%%%%%%%%%%%%%%%%%%%%%%%%%%%%%%%%%%%%%%%%%%%%%%%%%%%%%%%%%%%%%%%%%%%%%%%%%%%%%%%
%%%%%%%%%%%%%%%%%%%%%%%%%%%%%%%%%%%%%%%%%%%%%%%%%%%%%%%%%%%%%%%%%%%%%%%%%%%%%%%
% Abstract (250 words max)
%%%%%%%%%%%%%%%%%%%%%%%%%%%%%%%%%%%%%%%%%%%%%%%%%%%%%%%%%%%%%%%%%%%%%%%%%%%%%%%
%%%%%%%%%%%%%%%%%%%%%%%%%%%%%%%%%%%%%%%%%%%%%%%%%%%%%%%%%%%%%%%%%%%%%%%%%%%%%%%

% ****************************************************************************%
% To show Review Response highlighting and margin notes,
% \def\revnotes{}
% (see RevNotes.tex)
%\revOne{1}{Modified or added text} % To correspond to Reviewer #1 comment 1
%\revTwo{7}{Modified or added text} % To correspond to Reviewer #2 comment 7
% ****************************************************************************%

\begin{abstract}
We present variants of the Conjugate Gradient (CG), Conjugate Residual (CR), 
and Generalized Minimal Residual (GMRES) methods which are both \emph{pipelined} and \emph{flexible}.
These allow computation of inner products and norms to be overlapped with operator 
and nonlinear or nondeterministic preconditioner application.
The methods are hence aimed at hiding network latencies and synchronizations which can become
computational bottlenecks in Krylov methods on extreme-scale systems or in the 
strong-scaling limit. 
The new variants are not arithmetically equivalent to their base flexible Krylov methods,
but are chosen to be similarly performant in a realistic use case, the application 
of strong nonlinear preconditioners to large problems which require many Krylov iterations.
We provide scalable implementations of our methods as contributions to the \textsc{PETSc} package
and demonstrate their effectiveness with practical examples derived from 
models of mantle convection and lithospheric dynamics with heterogeneous
viscosity structure.
These represent challenging problems where
multiscale nonlinear preconditioners are required for the current state-of-the-art
algorithms, and are hence amenable to acceleration with our new techniques.
Large-scale tests are performed in the strong-scaling regime on a contemporary leadership supercomputer, 
where speedups approaching, and even exceeding $2\times$ can be observed.
We conclude by analyzing our new methods with a performance model targeted at future exascale machines.
\end{abstract}
\ifx\undefined\nosiam


\pagestyle{myheadings}
\thispagestyle{plain}
\markboth{P.~Sanan, S.M.~Schnepp, and D.A.~May}{Pipelined, Flexible Krylov Subspace Methods}
\fi

%%%%%%%%%%%%%%%%%%%%%%%%%%%%%%%%%%%%%%%%%%%%%%%%%%%%%%%%%%%%%%%%%%%%%%%%%%%%%%%
%%%%%%%%%%%%%%%%%%%%%%%%%%%%%%%%%%%%%%%%%%%%%%%%%%%%%%%%%%%%%%%%%%%%%%%%%%%%%%%
\section{Introduction}
%%%%%%%%%%%%%%%%%%%%%%%%%%%%%%%%%%%%%%%%%%%%%%%%%%%%%%%%%%%%%%%%%%%%%%%%%%%%%%%
%%%%%%%%%%%%%%%%%%%%%%%%%%%%%%%%%%%%%%%%%%%%%%%%%%%%%%%%%%%%%%%%%%%%%%%%%%%%%%

The current High Performance Computing (HPC) paradigm involves
computation on larger and larger clusters of individual compute nodes connected through a network.
Since single-core performance has plateaued, increased parallelism is used to increase total performance.
As systems scale in this way, new performance bottlenecks emerge in algorithms popularized before massive parallelism became relevant.
While peak compute power is abundant on modern high performance clusters,
utilizing them at full capacity requires modifications to many well established algorithms. 
We investigate modifications to flexible Krylov subspace methods to
better exploit computational resources available from modern clusters.  
As parallelism increases, collective operations are more likely to become computational bottlenecks.
For instance, distributed dot products can
become performance limiting due to network latency exposed by 
synchronized all-to-all communication.
Similar concerns apply to hypothetical future machines and to current machines in 
the strong-scaling limit, when small amounts of local work expose bottlenecks in 
collective operations.
One approach to mitigate these bottlenecks is to overlap these collective operations 
with local work, hiding the latency.
This is a non-trivial task due to data dependencies
which forbid arbitrary rearrangements of operations.

%=============================================================================%
\subsection{Pipelining Flexible Krylov Subspace Methods}
%=============================================================================%
A classical way to use a system with independent hardware resources at higher 
efficiency is to \emph{pipeline}.
Broadly, this technique can accelerate throughput of a repeated multi-stage process
when different stages require different resources.
Concurrent work is performed on multiple iterations of the process, 
overlapped such that multiple resources are used simultaneously.
This induces latency as the pipeline is filled as well as other overheads related 
to rearranging an algorithm and allowing simultaneous progress on multiple tasks.

We consider pipelining a particular class of algorithms, \revOne{11}{namely} Krylov 
subspace methods for the solution of large, sparse, linear systems $Ax=b$ \cite{Saad2003,VanderVorst2003}. 
Recent research has provided algorithms which loosen data dependencies in these algorithms, 
allowing communication- and computation-oriented resources to operate 
concurrently \cite{Ghysels2013,Ghysels2014,Hoemmen2010}.

Krylov subspace methods for scientific applications are typically preconditioned \cite{Wathen2015}.
This is essential for scalability in many cases.
For instance, when the operator $A$ is a discretized elliptic operator, a
multilevel preconditioner \cite{Smith2004} is required for algorithmic scalability.
Preconditioners are often not available as assembled matrices. 
They commonly involve nested approximate solves, an important case being the application of another Krylov method \cite{mszm2014}. 
The application of a nested Krylov method to a given tolerance and/or iteration count is
 a nonlinear operation; the optimal polynomial chosen at each iteration depends 
 on $A$, $b$, and $x_0$. 
 This requires the `outer' Krylov method to be \emph{flexible},
able to operate with a nonlinear preconditioner.
Flexible Krylov methods do not involve Krylov subspaces but approximations to them.
In comparison with their non-flexible counterparts, flexible Krylov methods typically have higher storage requirements,
as they cannot exploit the structure of true Krylov spaces.

In this work, we explore, analyze, and implement Krylov methods which are both pipelined and flexible. 
As shown in~\S\ref{sec:implementation}-\ref{sec:perf-model}, promising applications of the new
methods presented here involve \revOne{1}{preconditioners with nested inexact solves} used over many Krylov iterations. 
For the overlap of reductions to be useful, systems must involve large numbers of nodes, or one must be in a
strong-scaling regime where the time for a global reduction is comparable to the time spent performing
local and neighbor-wise work, or when local processing times are very variable.

Our methods are applicable to the solution of any system of linear equations,
but are expected to be useful when the number of processing units is large, relative 
to the local problem size, and when strong, complex, nonlinear (or even nondeterministic) 
preconditioners are applied to systems difficult enough to require enough Krylov 
iterations to amortize the overhead of filling a pipeline.
Thus, for the examples presented in \S\ref{sec:tests}, we focus on applications
in geophysics, in particular on mantle convection and lithospheric dynamics with
heterogeneous viscosity structure. 
These represent challenging problems, where strong
multiscale nonlinear preconditioners are typically required for acceptable convergence.

Nonlinear preconditioning is still an active area of research, and newly developed 
nonlinear methods \cite{Brune2015} may provide additional use cases for the methods presented here.
The advent of non-deterministic and randomized preconditioning techniques, promising 
for use with hybrid or heterogeneous clusters featuring accelerators, may provide future use cases for algorithms of the type presented here; these allow overlap of communication and computation and also loosen synchronization requirements, crucial when local processing times have heavy-tailed distributions \cite{Morgan2015}.

%=============================================================================%
\subsection{Notation and Algorithmic Presentation}\label{sec:notation}
%=============================================================================%
We use similar notation to analyze variants of the Conjugate Gradient (CG), Conjugate Residual (CR),
and Generalized \revOne{12}{Minimal} Residual (GMRES) methods, as described in Tables~\ref{tab:cgnotation} and
\ref{tab:gmresnotation}. 
This allows similar presentation of the methods, despite the fact that CG and CR methods are 
written as left preconditioned, while flexible GMRES methods are right preconditioned.
Operators are denoted by capital letters, vectors by lowercase Latin letters, and scalars by Greek letters.
Algorithms are presented using $0$-based iteration counts. 
For readability and analysis, main iteration loops are arranged such that variables indexed by the current index are updated. 
Convergence tests and optimizations related to initialization, elimination of
intermediates, and loop ordering are omitted for readability; the interested
reader is referred to \S\ref{sec:implementation} where full
open-source implementations are discussed.
  We adopt the following coloring scheme for operations involving vectors. 
  Scalar or `no operation' operations are not colored.
  %***************************************************************************%
  % NOTE : The text colors and symbols are hard-coded here (and nowhere else). 
  % If the colors are changed enough to have different names, they need to be 
  %  updated here as well, and similarly for the symbols.
  %***************************************************************************%
  \ifx\undefined\nosiam\begin{romannum}\else\begin{itemize}\fi
  \item \textLoc{Green} denotes an algorithmic step which is ``local'', 
  meaning that no internode communication is required in a typical parallel implementation.
  \item \textNei{Orange} denotes an algorithmic step which operates 
  ``neighbor-wise'', only involving communication between $O(1)$ neighboring nodes
   interspersed with local work. 
   Sparse matrix-vector multiplication typically falls into this category.
  \item \textRed{Red} denotes an algorithmic step which is ``global'',
   involving global communication, exposing the full latency of the network.
    These operations, which include computing vector dot products and norms, 
    can become the computational bottleneck in massively parallel systems.
  \item \textPC{Purple} is used for preconditioner application, which may 
  be local (e.g.~applying a Jacobi preconditioner), neighbor-wise (e.g.~applying 
  a sparse approximate inverse), or global (e.g.~applying a nested Krylov
  method).
  \ifx\undefined\nosiam\end{romannum}\else\end{itemize}\fi
  \revTwo{2}{Algorithm line numbers corresponding to local, neighbor-wise, global, and preconditioning operations
  are also marked with a dash(\textendash), plus(+), bullet(\textbullet),
  and asterisk(\textasteriskcentered), respectively.}

%=============================================================================%
\subsection{Contributions}
%=============================================================================%
\ifx\undefined\nosiam\begin{romannum}\else\begin{itemize}\fi
\item 
We present pipelined variants of flexible methods, allowing 
for inexact and variable preconditioning. (\S\ref{sec:pipel-flex-conj}-\ref{sec:gmres})
\item We provide open-source, scalable implementations of three solvers as contributions to the \textsc{PETSc} package. 
(\S\ref{sec:implementation})
\item We provide analysis of a novel modification of a ``naive'' pipelining which gives unacceptable convergence behavior in typical applications.
This leads to useful, pipelined variants of the Flexible CG (FCG) and Generalized CR (GCR) methods. (\S\ref{sec:pipefcg})
\item We demonstrate the use of the methods with
applications to solving linear systems from challenging problems in lithospheric
dynamics (\S\ref{sec:tests}), and with performance models to extrapolate to
exascale (\S\ref{sec:perf-model}).
\ifx\undefined\nosiam\end{romannum}\else\end{itemize}\fi

%%%%%%%%%%%%%%%%%%%%%%%%%%%%%%%%%%%%%%%%%%%%%%%%%%%%%%%%%%%%%%%%%%%%%%%%%%%%%%%
%%%%%%%%%%%%%%%%%%%%%%%%%%%%%%%%%%%%%%%%%%%%%%%%%%%%%%%%%%%%%%%%%%%%%%%%%%%%%%%
\section{Pipelined Flexible Conjugate Gradient Methods}
%%%%%%%%%%%%%%%%%%%%%%%%%%%%%%%%%%%%%%%%%%%%%%%%%%%%%%%%%%%%%%%%%%%%%%%%%%%%%%%
%%%%%%%%%%%%%%%%%%%%%%%%%%%%%%%%%%%%%%%%%%%%%%%%%%%%%%%%%%%%%%%%%%%%%%%%%%%%%%%
\label{sec:pipel-flex-conj}

\paragraph{Notation}\label{sec:cgnotation}
We base the notation used in this section, collected in Table \ref{tab:cgnotation}, 
on standard notation for Conjugate Gradient methods and its pipelined variant \cite{Ghysels2014}.
\begin{table}[h!]
\footnotesize %from SIAM latex guidelines
\begin{tabular}{l|l}
  $A$          & linear operator\\
  $M^{-1}$     & linear, symmetric positive definite, left preconditioner\\
  $B$          & nonlinear left preconditioner\\
  $x$          & true solution vector $A^{-1}b$\\
  $x_i$        & approximate solution vector at iteration $i = 0,1,\ldots$ \\ 
  $e_i$        & error $x - x_i$ at iteration $i$ \\
  $r_i$        & residual $b - A x_i = Ae_i$ at iteration $i$ \\
  $u_i$        & preconditioned residual $B(r_i)$ or $M^{-1}r_i$ \\
  $\tilde u_i$ & approximation to $u_i$, exact for a linear preconditioner\\ 
  $w_i$        & pipelining intermediate $Au_i$ or $A\tilde u_i$ \\
  $m_i$        & pipelining intermediate $B(w_i)$ or $M^{-1}w_i$\\
  $n_i$        & pipelining intermediate $Am_i$ \\
  $p_i$        & search direction or basis vector\\
  $s_i$        & transformed search direction $Ap_i$ \\
  $q_i$        & pipelining intermediate $B(s_i)$ or $M^{-1}s_i$ \\
  $z_i$        & pipelining intermediate $Aq_i$ \\
  $\alpha_i$ & scalar weight in solution and residual update\\
  $\beta_i$ & scalar weight in computation of new search direction or basis vector\\
  $\eta_i$ & squared norm $\langle Ap_i,p_i\rangle = ||p_i||_A^2$ \\
  $\delta_i$ & squared norm $\langle Au_i,u_i\rangle = ||u_i||_A^2$ \\
  $\gamma_i$ & $L_2$-inner product involving preconditioned or 
  unpreconditioned residuals\\ & depending on the method under consideration  
  ($\langle r_i,r_i\rangle$, $\langle u_i,u_i\rangle$, or $\langle u_i,r_i\rangle$)\\
\end{tabular}
\caption{Notation for Conjugate Gradient and related methods}
\label{tab:cgnotation}
\end{table}

%=============================================================================%
\subsection{Review of Conjugate Gradient Methods}
%=============================================================================%
\label{sec:revi-conj-grad}
\subsubsection{The Method of Conjugate Directions}\label{sec:meth-conj-direct}
A general class of methods known as 
\emph{conjugate direction (CD) methods} \cite{fox1948notes,Hestenes1952} have been investigated to compute (approximate) numerical 
solutions to the linear system $Ax=b$, where $A=A^H$ and $A$ is positive definite. 
Presume one has a set of $A$-orthogonal (``conjugate'') vectors 
$\{p_j\},j=0,\ldots,n-1$ available, that is a set of $n$ vectors with
the property $\langle p_j, p_k\rangle_A = \langle p_j ,A p_k \rangle = 0$ when $j \ne k$.
At each step of the algorithm, the conjugate directions method computes an approximate
 solution $x_i$ with the property that $||x - x_i||_A^2$ is 
 minimal for $x_i - x_0 \in \text{span}(p_0,\ldots,p_{i-1})$.
The resulting process is described in Algorithm \ref{alg:cd}.
Due to the unspecified nature of the $A$-orthogonal directions $\{p_j\}$,
it includes a wide range of algorithms, including Gaussian elimination 
and the family of Conjugate Gradient methods which are the main subject of this section.

We note an important property of all CD methods. By the minimization property, we have 
\begin{equation}\label{eq:cdorth}
e_i \perp_A \text{span}(p_0,\ldots,p_{i-1}) \iff r_i \perp \text{span}(p_0,\ldots,p_{i-1}).
\end{equation}
That is, the residual is orthogonal to the previously explored space in the standard norm
\footnote{Other norms may be used, as described by Hestenes \cite{Hestenes1956}.}.

The directions $p_j$ need not be specified before the algorithm begins; they can be computed as needed, based on the progress of the algorithm. 

\begin{algorithm}[!htbp]
\caption{Conjugate Directions \cite{Hestenes1952}}
\label{alg:cd}
\begin{algorithmic}[1]
\FunctionDef{CD}{$A$, $b$, $x_0$, $p_0,\ldots,p_{n-1}$} 
\StateNei{$r_0 \gets b - Ax_0$} %This is orange, but rarely actually costs a matmult, since either x_0 is 0, or you are restarting something, in which case you probably have r available
\StateLoc{$p_0 \gets r_0$}
\StateNei{$s_0 \gets Ap_0$}
\StateRed{$\gamma_0 \gets \langle p_0,r_0 \rangle$}
\StateRed{$\eta_0 \gets \langle p_0,s_0\rangle$}
\StateDef{$\alpha_0 \gets \gamma_0 / \eta_0$}
\ForDef{$ i = 1,2,\ldots,n$}
\StateLoc{$x_i \gets x_{i-1} + \alpha_{i-1}p_{i-1}$}
\StateLoc{$r_i \gets r_{i-1} - \alpha_{i-1}s_{i-1}$}
\StateRed{$\gamma_i \gets \langle p_i,r_i \rangle $}
\StateNei{$s_i \gets Ap_i$}
\StateRed{$\eta_i \gets \langle s_i,p_i\rangle$}
\StateDef{$\alpha_i \gets \gamma_i / \eta_i$}
\EndFor
\EndFunction
\end{algorithmic}
\revOne{13}{}
\end{algorithm}

\subsubsection{Flexible Conjugate Gradient Methods}\label{sec:flex-conj-grad}

To ensure rapid convergence of the CD algorithm, it is desirable at each iteration 
for the new search direction $p_i$ to be well-aligned with the remaining error $e_i = x-x_i $. 
As the true error is unknown, one available option is the residual $r_i = Ae_i$,
 $A$-orthogonalized against the other search directions.

The residual can be seen as a transformed error or as a gradient descent direction
 for the function $\tfrac{1}{2}||x||^2_A - \langle b,x\rangle$, a minimizer of which 
 is a solution \revOne{14}{of} $Ax=b$.
 Preconditioning can be seen as an attempt to obtain
 a direction better aligned with $e_i$ by instead
 $A$-orthogonalizing $u_i \doteq B(r_i) = B(Ae_i)$ 
 against the previous search directions $p_0,\ldots,p_{i-1}$. 
If $B$ is an approximate inverse of $A$, it is effective by this criterion,
and many preconditioners are motivated as such.

The process just described, found in Algorithm \ref{alg:fcg}, describes the 
Flexible Conjugate Gradient (FCG) method \cite{Notay2000}, with complete orthogonalization, 
i.e., the new search direction is found by $A$-orthogonalizing the residual 
against all previous directions in every iteration. 
Full $A$-orthogonalization requires potentially excessive memory usage and computation. 
Hence, it is a common approach to only $A$-orthogonalize against a number $\nu_i$ of previous directions. 
Note however that only for full orthogonalization ($\nu_i = i$) is a true
$A$-orthogonal basis computed by the Gram-Schmidt process in 
lines~\ref{algline:GSstart}--\ref{algline:GSend} of Algorithm~\ref{alg:fcg}, 
defining a CD method in the strict sense.

%Figure \ref{fig:fcg} shows a schematic of the data dependencies in the algorithm. 
%In this and similar figures, arrows are colored based on the type of operation they depend on. 
%Subscripts have been suppressed, and braces surround quantities for which values from previous iterations are required.

One can observe from Algorithm~\ref{alg:fcg} that no overlapping of reductions 
and operator or preconditioner application is possible;
the dot products depend on the immediately preceding computations.

\subsubsection{Preconditioned Conjugate Gradients} \label{sec:pcg}
Let $\Kry^j(A,b)$ represent the $j$th Krylov subspace $\text{span}(b,Ab,\ldots,A^{j-1}b)$.
If $B$ is a linear operator $B(v) \equiv
M^{-1}v$, then line~\ref{algline:fcgorthog} of Algorithm \ref{alg:fcg} evaluates to zero for all but the most recent previous direction, and the search directions lie in Krylov subspaces
$\Kry^i(M^{-1}A,M^{-1}b)$.
We thus recover the Preconditioned Conjugate Gradient Method \cite{Golub2007,Hestenes1956,Shewchuk1994} shown in Algorithm \ref{alg:pcg}. 
%The main loop is represented graphically in Figure~\ref{fig:pcg}. 
Again, communication and computation cannot be overlapped using the algorithm as written, and in a parallel implementation, two reductions must be performed per iteration.

\ifx\undefined\nosiam\ifx\undefined\nosiam\capstartfalse\fi\fi % to get hypcap to ignore this wrapper figure*
\begin{figure*} % wrap in figure* so this floats
  \begin{minipage}[t]{0.48\textwidth}
    \begin{algorithm}[H]
      \caption{Flexible Conjugate Gradients \cite{Notay2000} \white{xxxxx}} %invis x's to align titles
      \label{alg:fcg}
      \begin{algorithmic}[1]
        \FunctionDef{FCG}{$A$, $B$, $b$, $x_0$} 
        \StateNei{$r_0 \gets b - Ax_0$}
        \StatePC{$u_0 \gets B(r_0)$}
        \StateLoc{$p_0 \gets u_0$}
        \StateNei{$s_0 \gets Ap_0$}
        \StateRed{$\gamma_0 \gets \langle u_0,r_0 \rangle$}
        \StateRed{$\eta_0 \gets \langle p_0,s_0 \rangle$}
        \StateDef{$\alpha_0 \gets \gamma_0 / \eta_0$}
        \ForDef{$ i = 1,2,\ldots$}
        \StateLoc{$x_i \gets x_{i-1} + \alpha_{i-1}p_{i-1}$}
        \StateLoc{$r_i \gets r_{i-1} - \alpha_{i-1}s_{i-1}$}
        \StatePC{$u_i \gets B(r_i)$}
        \StateRed{$\gamma_i \gets \langle u_i,r_i \rangle $}
        \ForDef {$k = i-\nu_i,\ldots,i-1$} \label{algline:GSstart}
        \StateRed{$\beta_{i,k} \gets \tfrac{-1}{\eta_k}\langle u_i,s_k\rangle$} \label{algline:fcgorthog}
        \EndFor
        \StateLoc{$p_i \gets u_i + \sum_{k=i-\nu_i}^{i-1}\beta_{i,k}p_k$} \label{algline:GSend}
        \StateNei{$s_i \gets Ap_i$}
        \StateRed{$\eta_i \gets \langle p_i,s_i\rangle$}
        \StateDef{$\alpha_i \gets \gamma_i / \eta_i$}
        \EndFor
        \EndFunction
      \end{algorithmic}
    \end{algorithm}
  \end{minipage}%
  \hfill
  \begin{minipage}[t]{0.48\textwidth}
    \vspace{0pt}  
    \begin{algorithm}[H]
      \caption{Preconditioned Conjugate Gradients \cite{Hestenes1956} \label{alg:pcg}}
      \begin{algorithmic}[1]
        \FunctionDef{PCG}{$A$, $M^{-1}$, $b$, $x_0$} 
        \StateNei{$r_0 \gets b - Ax_0$} %This is orange, but rarely actually costs a matmult, since either x_0 is 0, or you are restarting something, in which case you probably have r
        \StatePC{$u_0 \gets M^{-1}r_0$}
        \StateLoc{$p_0 \gets u_0$}
        \StateNei{$s_0 \gets Ap_0$}
        \StateRed{$\gamma_0 \gets \langle u_0,r_0 \rangle$}
        \StateRed{$\eta_0 \gets \langle s_0,p_0 \rangle$}
        \StateDef{$\alpha_0 \gets \gamma_0 / \eta_0$}
        \ForDef{$ i = 1,2,\ldots$}
        \StateLoc{$x_i \gets x_{i-1} + \alpha_{i-1}p_{i-1}$}
        \StateLoc{$r_i \gets r_{i-1} - \alpha_{i-1}s_{i-1}$}
        \StatePC{$u_i \gets M^{-1}r_i$}
        \StateRed{$\gamma_i \gets \langle u_i,r_i \rangle $}
        \StateDef{} %empty for alignment
        \StateDef{$\beta_i \gets \gamma_i / \gamma_{i-1}$}
          \vspace{4pt} %for alignment
        \StateLoc{$p_i \gets u_i + \beta_ip_{i-1}$}
        \StateNei{$s_i \gets Ap_i$}
        \StateRed{$\eta_i \gets \langle s_i,p_i\rangle$}
        \StateDef{$\alpha_i \gets \gamma_i / \eta_i$}
        \EndFor
        \EndFunction
      \end{algorithmic}
    \end{algorithm}
  \end{minipage}%
  \revTwo{9}{}
\end{figure*}
\ifx\undefined\nosiam\capstarttrue\fi % to get hypcap to ignore this wrapper figure*

\subsubsection{Pipelined Preconditioned Conjugate Gradients}

Any rearrangement of the PCG algorithm beyond trivial options
must appeal to a notion of algorithmic equivalence more general than producing identical floating-point results. 
These notions include algebraic rearrangements, which typically do not produce 
the same iterates in finite precision arithmetic, but which are equivalent in exact arithmetic. 
One such operation is to rearrange the CG algorithm to involve a single
reduction, requiring less communication and synchronization overhead on parallel systems.
This was done by Chronopoulos and Gear \cite{Chronopoulos1989}, as described in Algorithm \ref{alg:cgcg}.
 The rearrangement is also particularly important in the context of implementing efficient Krylov methods for current hybrid
systems involving high throughput coprocessors for which kernel fusion is beneficial \cite{Rupp2014}.

This rearrangement is possible because CG relies on the geometric structure induced by interpreting symmetric positive definite operators as inner products.
Noting that $\eta_i$ is defined as the squared $A$-norm of $p_i$ and that $p_i$ is constructed by $A-$orthogonalizing $u_i$ with respect to $p_{i-1}$, one can use the Pythagorean Theorem to write
\begin{equation*}
\eta_i = ||p_i||_A^2 = ||u_i + \beta p_{i-1}||_A^2 = ||u_i||_A^2 - |\beta_i|^2||p_{i-1}||_A^2 = ||u_i||_A^2 - |\beta_i|^2 \eta_{i-1}
\end{equation*}
where terms involving products of directions at different iteration numbers
are zero because of the $A$-orthogonality of $p_i$ and $p_{i-1}$. 
 The quantity $\delta_i \doteq ||u_i||_A^2 = \langle u_i,Au_i\rangle $  does not depend on $p_i$ and can thus be computed at the same time as $\gamma_i$. The quantity $w_i \doteq Au_i$ needs to be computed, which apparently adds another matrix multiply to the overall iteration. However, a second type of structure, namely linear structure, can be exploited to reuse this computation and remove the existing matrix multiply with a recurrence relation
\begin{equation}\label{eq:linearrecursion}
s_i = Ap_i = Au_i + \beta_i Ap_{i-1} = w_i + \beta s_{i-1}.
\end{equation}
These rearrangements require additional storage (one extra vector $w$
in this case) and floating point computations and come with potential
numerical penalties as the computed norm could become negative (``\revOne{2}{norm}
breakdown'') in finite precision arithmetic. 
These breakdowns are most easily handled through algorithm restarts. 
Also, errors can accumulate in recursively computed variables.

As $\eta_i$ is only used to calculate $\alpha_i$, one can avoid computing it and update $\alpha_i$ directly as
\begin{equation*}
\alpha_i = \frac{\gamma_i}{\eta_i} = \frac{\gamma_i}{(\delta_i - |\beta_i|^2 \eta_{i-1})} = \frac{\gamma_i}{\delta_i - \left(\frac{\gamma_{i}}{\gamma_{i-1}}\right)^2 \left( \frac{\gamma_{i-1}}{\alpha_{i-1}}\right) } =  \frac{\gamma_i}{\delta_i - \bar\beta_i \frac{\bar\gamma_i}{\alpha_{i-1}}}.
\end{equation*}
This recurrence is more concise, but we retain explicit computation of $\eta_i$
in the algorithms as presented here to allow for an easier comparison.
\ifx\undefined\nosiam\capstartfalse\fi % to get hypcap to ignore this wrapper figure*
\begin{figure*} % wrap in figure* so this floats
\begin{minipage}{0.48\textwidth}
\vspace{0pt}  
\begin{algorithm}[H]
  \caption{Chronopoulos-Gear Conjugate Gradients \cite{Chronopoulos1989}}
  \label{alg:cgcg}
    \begin{algorithmic}[1]
      \FunctionDef{CGCG}{$A$, $M^{-1}$, $b$, $x_0$} 
      \StateNei{$r_0 \gets b-Ax_0$}
      \StatePC{$u_0 \gets M^{-1}r_0$}
      \StateNei{$w_0 \gets Au_0$}
      \StateDef{} % empty, for alignment
      \StateDef{} % empty, for alignment
      \StateRed{$\gamma_0 \gets \langle u_0,r_0 \rangle $}
      \StateRed{$\delta_0 \gets \langle u_0,w_0 \rangle $}
      \StateDef{$\eta_0 \gets \delta_0 $}
      \StateLoc{$p_0 \gets u_0$}
      \StateLoc{$s_0 \gets w_0$}
      \StateDef{} % empty, for alignment
      \StateDef{} % empty, for alignment
      \StateDef{$\alpha_0 \gets \gamma_0 / \eta_0$}
      \ForDef{$ i = 1,2,\ldots$}
      \StateLoc{$x_i \gets x_{i-1} + \alpha_{i-1}p_{i-1}$}
      \StateLoc{$r_i \gets r_{i-1} - \alpha_{i-1}s_{i-1}$}
      \StateDef{} % empty, for alignment
      \StateDef{} % empty, for alignment
      \StatePC{$u_i \gets M^{-1}r_i$} \label{algline:cgcgu}
      \StateNei{$w_i \gets Au_i$}
      \StateRed{$\delta_i \gets \langle u_i,w_i \rangle $}
      \StateRed{$\gamma_i \gets \langle u_i,r_i \rangle $}
      \StateDef{$\beta_i \gets \gamma_i  / \gamma_{i-1}$}
      \StateDef{$\eta_i \gets \delta_i - |\beta_i|^2\eta_{i-1}$}
      \StateDef{$\alpha_i \gets \gamma_i / \eta_i$}
      \StateLoc{$p_i \gets u_i + \beta_{i}p_{i-1}$}
      \StateLoc{$s_i \gets w_i + \beta_is_{i-1}$}
      \StateDef{} % empty, for alignment
      \StateDef{} % empty, for alignment
      \EndFor
      \EndFunction
      \vspace{3pt}
      %\ifx\undefined\nosiam\vspace*{95pt}\else\vspace*{85pt}\fi
    \end{algorithmic}
  \end{algorithm}
\end{minipage}%
\vspace{0pt}  
\hspace{3px}
\begin{minipage}{0.48\textwidth}
\begin{algorithm}[H]
  \caption{Pipelined Conjugate Gradients \cite{Ghysels2014} \white{xxxxx}}
  \label{alg:pipecg}
  \begin{algorithmic}[1]
    \FunctionDef{PIPECG}{$A$, $M^{-1}$, $b$, $x_0$} 
    \StateNei{$r_0 \gets b - Ax$}
    \StatePC{$u_0 \gets M^{-1}r_0$}
    \StateNei{$w_0 \gets Au_0$}
    \StatePC{$m_0 \gets M^{-1}w_0$}
    \StateNei{$n_0 \gets Am_0$}
    \StateRed{$\gamma_0 \gets \langle u_0,r_0 \rangle $}
    \StateRed{$\delta_0 \gets \langle u_0,w_0 \rangle $}
    \StateDef{$\eta_0 \gets \delta_0$}
    \StateLoc{$p_0 \gets u_0$}
    \StateLoc{$s_0 \gets w_0$}
    \StateLoc{$q_0 \gets m_0$}
    \StateLoc{$z_0 \gets n_0$}
    \StateDef{$\alpha_0 \gets \gamma_0 / \eta_0$}
    \ForDef{$ i = 1,2,\ldots$}
    \StateLoc{$x_i \gets x_{i-1} + \alpha_{i-1}p_{i-1}$}
    \StateLoc{$r_i \gets r_{i-1} - \alpha_{i-1}s_{i-1}$}
    \StateLoc{$u_i \gets u_{i-1} - \alpha_{i-1}q_{i-1}$}
    \StateLoc{$w_i \gets w_{i-1} - \alpha_{i-1}z_{i-1}$}
    \StatePC{$m_i \gets M^{-1}w_i$}
    \StateNei{$n_i \gets Am_i$}
    \StateRed{$\gamma_i \gets \langle u_i,r_i \rangle $}
    \StateRed{$\delta_i \gets \langle u_i,w_i \rangle $}
    \StateDef{$\beta_i \gets \gamma_i  / \gamma_{i-1}$}
    \StateDef{$\eta_i \gets \delta_i - |\beta_i|^2\eta_{i-1}$}
    \StateDef{$\alpha_i \gets \gamma_i / \eta_i$}
    \StateLoc{$p_i \gets u_i + \beta_{i}p_{i-1}$}
    \StateLoc{$s_i \gets w_i + \beta_is_{i-1}$}
    \StateLoc{$q_i \gets m_i + \beta_i q_{i-1}$}
    \StateLoc{$z_i \gets n_i + \beta_i z_{i-1}$}
    \EndFor
    \EndFunction
  \end{algorithmic}
  \end{algorithm}
\end{minipage}%
\revTwo{9}{}
\end{figure*}
\ifx\undefined\nosiam\capstarttrue\fi % to get hypcap to ignore this wrapper figure*

The Chronopoulos-Gear CG method cannot overlap the preconditioner and sparse matrix 
multiply with the reductions, so some further rearrangement is desirable to allow 
concurrent use of computational resources. 
Ghysels and Vanroose \cite{Ghysels2014} developed a further variant of the Chronopoulos-Gear
 algorithm to accomplish this, as described in Algorithm~\ref{alg:pipecg}. 
 The algorithm takes advantage of multiple ``unrollings'' using the linear
 structure of the variable updates, as in \eqref{eq:linearrecursion}. 
For example, since the computation of $u_i = M^{-1}r_i$ on line \ref{algline:cgcgu} 
of Algorithm \ref{alg:cgcg} blocks the subsequent inner product involving $u_i$, one can use the identity
\begin{equation}\label{eq:pcunroll}
u_i = M^{-1}r_i = M^{-1}r_{i-1} - \alpha_{i-1}M^{-1} s_i = u_{i-1} - \alpha_{i-1}q_i
\end{equation}
which shifts the application of $M^{-1}$ from $r_i$ to $s_i$, to compute a new 
variable $q_i \doteq M^{-1}s_i$. Similarly, $q_i$ is rewritten until 
computation, for $m_i$ in Algorithm \ref{alg:pipecg}, can be 
overlapped with the dot products in the previous iteration. 
The identity in \eqref{eq:pcunroll} relies on the linearity of $M^{-1}$.
This rearrangement allows for more concurrent
use of computational resources at the cost of storing more vectors,
performing more floating point operations, and reducing numerical
stability \cite{Ghysels2014}. 
Examples as well as a performance model 
\cite{Ghysels2013,Ghysels2014} show performance gains on current and future
parallel systems. 

The Pipelined Conjugate Gradient method as presented in
Algorithm~\ref{alg:pipecg} requires storing 10 vectors
compared to 6 for PCG and 7 for CGCG.

%=============================================================================%
\subsection{Pipelined Flexible Conjugate Gradient Methods} \label{sec:pipefcg}
%=============================================================================%

The FCG algorithm can be modified in the same manner used to produce the Chronopoulos-Gear CG algorithm. This does not involve any assumption of linearity of $B$, so the new single reduction FCG algorithm described in Algorithm \ref{alg:cgfcg} is equivalent to FCG in exact arithmetic, for an arbitrary preconditioner.

% CGFCG
\begin{algorithm}[!htbp]
\caption{Single Reduction Flexible Conjugate Gradients}
\label{alg:cgfcg}
\begin{algorithmic}[1]
\FunctionDef{CGFCG}{$A$, $B$, $b$, $x_0$} %TODO : revisit these algorithm names
\StateNei{$r_0 \gets b - Ax_0$} %This is orange, but rarely actually costs a matmult, since either x_0 is 0, or you are restarting something, in which case you probably have r
\StatePC{$u_0 \gets B(r_0)$}
\StateNei{$w_0 \gets Ap_0$}
\StateRed{$\gamma_0 \gets \langle u_0,r_0 \rangle$}
\StateRed{$\delta_0 \gets \langle u_0,w_0\rangle$}
\StateLoc{$p_0 \gets u_0$}
\StateLoc{$s_0 \gets w_0$}
\StateDef{$\eta_0 \gets \delta_0$}
\StateDef{$\alpha_0 \gets \gamma_0 / \eta_0$}
\ForDef{$ i = 1,2,\ldots$}
\StateLoc{$x_i \gets x_{i-1} + \alpha_{i-1}p_{i-1}$}
\StateLoc{$r_i \gets r_{i-1} - \alpha_{i-1}s_{i-1}$}
\StatePC{$u_i \gets B(r_i)$}
\StateNei{$w_i \gets Au_i$}
\StateRed{$\gamma_i \gets \langle u_i,r_i \rangle $}
\ForDef {$k = i-\nu_i,\ldots,i-1$}
\StateDef{\textRed{$\beta_{i,k} \gets \tfrac{-1}{\eta_k}\langle u_i,s_k\rangle$}}
\EndFor
\StateRed{$\delta_i \gets \langle u_i, w_i \rangle$}
\StateLoc{$p_i \gets u_i + \sum_{k=i-\nu_i}^{i-1}\beta_{i,k}p_k$}
\StateLoc{$s_i \gets w_i + \sum_{k=i-\nu_i}^{i-1}\beta_{i,k}s_k$}
\StateDef{$\eta_i \gets \delta_i - \sum_{k=i-\nu_i}^{i-1}\beta_{i,k}^2\eta_k$}
\StateDef{$\alpha_i \gets \gamma_i / \eta_i$}
\EndFor
\EndFunction
\end{algorithmic}
\end{algorithm}

\ifx\undefined\nosiam\capstartfalse\fi % to get hypcap to ignore this wrapper figure*
\begin{figure*}
\revTwo{3}{}
\begin{minipage}[t]{0.48\textwidth}
  \vspace{0pt}  
  \begin{algorithm}[H]
    \caption{Pipelined FCG (Naive)}
    \label{alg:naivepipefcg}
    \begin{algorithmic}[1]
      \FunctionDef{PIPEFCG?}{$A$, $B$, $b$, $x_0$} 
      \StateNei{$r_0 \gets b - Ax_0$} %This is orange, but rarely actually costs a matmult, since either x_0 is 0, or you are restarting something, in which case you probably have r
      \StatePC{$\tilde u_0 \gets B(r_0)$} %this one doesn't actually need a tilde
      \StateLoc{$p_0 \gets \tilde u_0$}
      \StateNei{$w_0 \gets Ap_0$} 
      \StatePC{$m_0 \gets B(w_0)$}
      \StateNei{$n_0 \gets Am_0$}
      \StateRed{$\gamma_0 \gets \langle \tilde u_0,r_0 \rangle$}
      \StateRed{$\delta_0 \gets \langle \tilde u_0, w_0\rangle$}
      \StateLoc{$s_0 \gets w_0$}
      \StateLoc{$q_0 \gets m_0$}
      \StateLoc{$z_0 \gets n_0$}
      \StateDef{$\eta_0 \gets \delta_0$}
      \StateDef{$\alpha_0 \gets \gamma_0 / \eta_0$}
      \ForDef{$ i = 1,2,\ldots$}
      \StateLoc{$x_i \gets x_{i-1} + \alpha_{i-1}p_{i-1}$}
      \StateLoc{$r_i \gets r_{i-1} - \alpha_{i-1}s_{i-1}$}
      \StateLoc{$\tilde u_i \gets \tilde u_{i-1} - \alpha_{i-1}q_{i-1}$} \label{algline:naivepipefcgu}
      \StateLoc{$w_i \gets w_{i-1} - \alpha_{i-1}z_{i-1}$} \label{algline:naivepipefcgw}
      \StateRed{$\gamma_i \gets \langle \tilde u_i,r_i \rangle $}
      \ForDef {$k=i--\nu_i,\ldots,i-1$}
      \StateRed{$\beta_{i,k} \gets \tfrac{-1}{\eta_k}\langle \tilde u_i,s_k\rangle$}
      \EndFor
      \StateRed{$\delta_i \gets \langle \tilde u_i, w_i \rangle$}
      \StatePC{$\mathbf{m_i \gets B(w_i)}$} \label{algline:naivepipefcgm}
      \StateNei{$n_i \gets Am_i$}
      \StateLoc{$p_i \gets \tilde u_i + \sum_{k=i-\nu_i}^{i-1}\beta_{i,k}p_k$}
      \StateLoc{$s_i \gets w_i + \sum_{k=i-\nu_i}^{i-1}\beta_{i,k}s_k$}
      \StateLoc{$q_i \gets m_i + \sum_{k=i-\nu_i}^{i-1}\beta_{i,k}q_k$}
      \StateLoc{$z_i \gets n_i +  \sum_{k=i-\nu_i}^{i-1}\beta_{i,k}z_k$}
      \StateDef{$\eta_i \gets \delta_i - \sum_{k=i-\nu_i}^{i-1}\beta_{i,k}^2\eta_k$}
      \StateDef{$\alpha_i \gets \gamma_i / \eta_i$}
      \EndFor
      \EndFunction
    \end{algorithmic}
  \end{algorithm}
\end{minipage}%
\hspace{3px}
\begin{minipage}[t]{0.48\textwidth}
  \vspace{0pt}
  \begin{algorithm}[H]
    \caption{Pipelined FCG \white{()}}
    \label{alg:modpipefcg}
    \begin{algorithmic}[1]
      \FunctionDef{PIPEFCG}{$A$, $B$, $b$, $x_0$} 
      \StateNei{$r_0 \gets b - Ax_0$} 
      \StatePC{$\tilde u_0 \gets B(r_0)$} %this one doesn't actually need a tilde
      \StateLoc{$p_0 \gets \tilde u_0$}
      \StateNei{$w_0 \gets Ap_0$}
      \StatePC{$m_0 \gets B(w_0)$}
      \StateNei{$n_0 \gets Am_0$}
      \StateRed{$\gamma_0 \gets \langle u_0,r_0 \rangle$}
      \StateRed{$\delta_0 \gets \langle u_0, w_0\rangle$}
      \StateLoc{$s_0 \gets w_0$}
      \StateLoc{$q_0 \gets m_0$}
      \StateLoc{$z_0 \gets n_0$}
      \StateDef{$\eta_0 \gets \delta_0$}
      \StateDef{$\alpha_0 \gets \gamma_0 / \eta_0$}
      \ForDef{$ i = 1,2,\ldots$}
      \StateLoc{$x_i \gets x_{i-1} + \alpha_{i-1}p_{i-1}$}
      \StateLoc{$r_i \gets r_{i-1} - \alpha_{i-1}s_{i-1}$}
      \StateLoc{$\tilde u_i \gets \tilde u_{i-1} - \alpha_{i-1}q_{i-1}$}
      \StateLoc{$w_i \gets w_{i-1} - \alpha_{i-1}z_{i-1}$}
      \StateRed{$\gamma_i \gets \langle \tilde u_i,r_i \rangle $}
      \ForDef {$k=i-\nu_i,\ldots,i-1$}
      \StateRed{$\beta_{i,k} \gets \tfrac{-1}{\eta_k}\langle \tilde u_i,s_k\rangle$}
      \EndFor
      \StateRed{$\delta_i \gets \langle \tilde u_i, w_i \rangle$}
      \StatePC{$\mathbf{m_i \gets \tilde u_i + B(w_i-r_i)}$} \label{algline:stabpipefcgm}
      \StateNei{$n_i \gets Am_i$}
      \StateLoc{$p_i \gets \tilde u_i + \sum_{k=i-\nu_i}^{i-1}\beta_{i,k}p_k$}
      \StateLoc{$s_i \gets w_i + \sum_{k=i-\nu_i}^{i-1}\beta_{i,k}s_k$}
      \StateLoc{$q_i \gets m_i + \sum_{k=i-\nu_i}^{i-1}\beta_{i,k}q_k$}
      \StateLoc{$z_i \gets n_i +  \sum_{k=i-\nu_i}^{i-1}\beta_{i,k}z_k$}
      \StateDef{$\eta_i \gets \delta_i - \sum_{k=i-\nu_i}^{i-1}\beta_{i,k}^2\eta_k$}
      \StateDef{$\alpha_i \gets \gamma_i / \eta_i$}
      \EndFor
      \EndFunction
    \end{algorithmic}
  \end{algorithm}
\end{minipage}
\end{figure*}
\ifx\undefined\nosiam\capstarttrue\fi % to get hypcap to ignore this wrapper figure*

One can now attempt to pipeline FCG in the same manner as pipelined variants of CG are obtained - 
this leads to the variant described in Algorithm \ref{alg:naivepipefcg}.
The quantity $\tilde u$ is used to reinforce the idea that the preconditioned search directions differ 
in general from those \revOne{15}{computed by FCG}.

A variant analogous to Gropp's asynchronous CG \cite{Gropp2010pres} could also be defined.
We do so for the corresponding method based on \revOne{16}{conjugate} residuals in \S\ref{sec:pipelined-GCR}.

\subsubsection{Issues with Naive Pipelining and Variable Preconditioning}
\label{sec:fcg-issues}

\begin{figure}
\begin{center}
\includegraphics[width=0.8\textwidth]{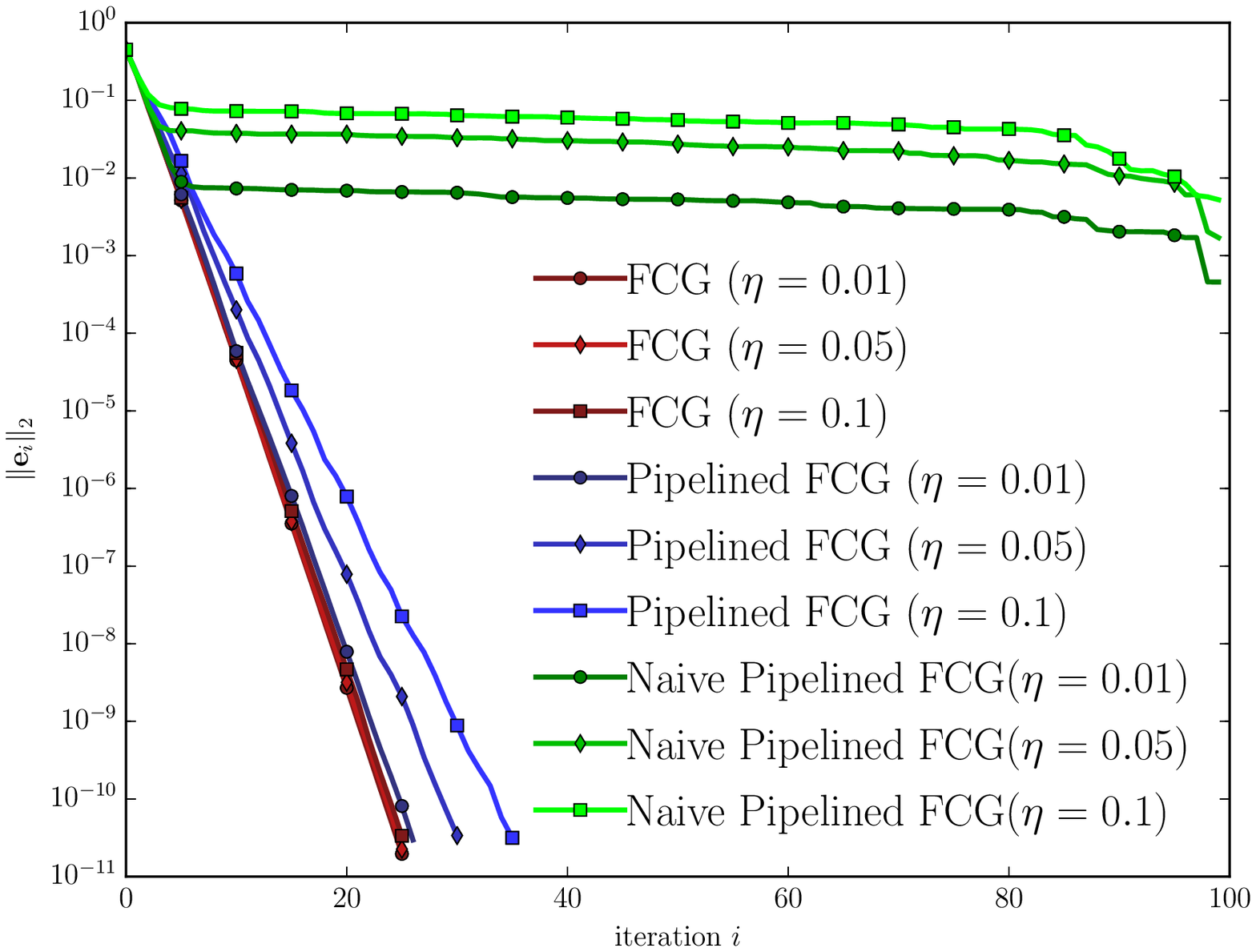}
\caption{A toy diagonal system, $n=100$, with condition number
 $5$, equally-spaced eigenvalues, and $b = \mathbb{1}/\sqrt{n}$,   
 solved with FCG, 
 naively pipelined FCG (Algorithm \ref{alg:naivepipefcg}), and our PIPEFCG (Algorithm \ref{alg:modpipefcg}).
 The pathological preconditioner $B$ adds Gaussian noise of magnitude $\eta||r||$ to the residual $r$.
  $\numax = n$. 
  Note that convergence for the naively pipelined algorithm stagnates at 
  a relative error of approximately $\eta$.}
\label{fig:FCGtoy}
\end{center}
\revOne{3}{}
\end{figure}

Convergence of Algorithm \ref{alg:naivepipefcg} typically
stagnates for nonlinear preconditioners $B$. 
In cases where the preconditioner
has an associated ``noise level'' or tolerance which characterizes its variability,
stagnation is typically observed at a relative error comparable to this quantity. 
This is illustrated in Figure \ref{fig:FCGtoy}, plotting convergence for a 
diagonal system with a preconditioner which simply adds Gaussian noise to the residual.

Analysis of preconditioned conjugate gradient methods typically makes use of the 
preconditioner as an inner product. 
The CG method can be seen as a Ritz-Galerkin approach to a 
Krylov subspace method where approximations in successive Krylov subspaces are 
chosen to satisfy the condition that the resulting residual is orthogonal to 
the Krylov subspace \cite[\S4.1]{VanderVorst2003}. 
Choosing a different definition of orthogonality, based on the metric induced 
by $M^{-1}$, produces the PCG method.
This is equivalent to the notion of simply replacing the $L_2$ inner product 
with a different inner product, wherever it appears in the algorithm.
This freedom to choose norms was recognized soon after the development of the 
CG method \cite{Hestenes1956}. In terms of the PCG, Algorithm \ref{alg:pcg},
\begin{equation}
  \label{eq:PCGinnerproduct}
  \gamma_i = \langle r_i,r_i \rangle_{M^{-1}} = \langle u_i, r_i\rangle.
  \end{equation}

Here, the nonlinearity of the preconditioner precludes direct use of this idea. 
We opt to consider $B$ in a more general sense of a preconditioner, as an 
approximate inverse to $A$. 
%PDS : This could be made more specific.. we don't need an approximate inverse, really. What we want is something like a scaled approximate inverse on a small number of independent  eigenspaces, approximately, so we have clustered eigs.
Specifically, noting that $Ae_i = r_i$, we shall later 
require that $||B(r_i) - e_i||/||e_i||$ is controlled.

Naively replacing the linear preconditioner $M^{-1}$ in~\eqref{eq:pcunroll} with a 
nonlinear one, we obtain the unrolling used in Algorithm \ref{alg:naivepipefcg},
\begin{equation}\label{eq:pcunroll2}
  u_i = B(r_i) \approx B(r_{i-1}) - B(\alpha_{i-1} s_i) \approx u_{i-1} - \alpha_{i-1}q_i,
\end{equation}
which is inexact even in exact arithmetic. 
This unrolling is not a priori unacceptable, as FCG is designed to cope with 
varying and inexact preconditioning, and the rearrangement amounts to a different
 nonlinear preconditioner, albeit one derived using intuition from the linear case.
However, unlike the unrolling~\eqref{eq:pcunroll} the above leads to an
accumulation of errors
in the preconditioned residual $\tilde u_i$, with respect to $u_i = B(r_i)$. 
This invalidates the equivalence~\eqref{eq:PCGinnerproduct}, 
which PCG structurally relies on for projecting out search directions from the residual. 

Note that when introducing a nonlinear preconditioner with FCG\footnote{We assume $\nu_i = i$, that is complete orthogonalization, 
for all of the analysis in this section.}
the basis vectors $p_j$ remain $A$-orthogonal with this rearrangement, yet the residuals $r_i$, 
which are $M^{-1}$-orthogonal when $B(v) \equiv M^{-1}v$, are now orthogonal with respect 
to an equivalent linear operator $\hat B$, using Notay's concept and notation of 
an \emph{equivalent linear operator} \cite{Notay2000}.
This operator is defined by the relations $u_i = B(r_i) = \hat B r_i$, depends on $b$ and $x_0$, and is 
unique except in degenerate cases.
%It should be noted, however, that $\hat B$ is not by definition symmetric or positive definite, so 
%does not define an inner product.
It seems to be the case that the equivalent linear operator $\hat{\tilde B}$ implicitly defined 
(Cf. \eqref{eq:pcunroll}) by the update on line~\ref{algline:naivepipefcgu} 
of Algorithm \ref{alg:naivepipefcg} is unbounded as $i \to \infty$,\revTwo{5}{}
mapping residuals of arbitrarily small magnitude to preconditioned residuals of the same order,
making it an unacceptable approximation for any given linear 
operator $M^{-1}$ bounded as $n \to \infty$. 

Thus, it will be useful to distinguish between the equivalent linear operators 
as induced by the base and pipelined methods, $\hat B$ and $\hat{\tilde B}$, respectively. 
The operator $\hat{\tilde B}$ is defined by the relations $\tilde u_i = \hat{\tilde B}r_i$.
The quantities $\alpha_{i}$, in a CD method like FCG, are chosen to minimize $||x- (x_{i+1} + \alpha_{i-1} u_i)||_A$. However, in the pipelined variant,
the effective update is $x_i \gets x_{i-1} + \alpha_{i-1} \tilde u_i$.
Thus, deviation from the base method induced by pipelining can of course 
be characterized by the difference between $u$ and $\tilde u$, 
or equivalently the difference between $\tilde B$ and $\hat{\tilde B}$. 
The advantage of this second view is that it suggests a repair procedure.
Namely, one should respect the already-defined part of $\hat B$ as much as possible when computing $\tilde u$, 
hence partially defining $\hat{\tilde B}$.
One approach is to note that
the operation of the equivalent linear operator in direction $r_i$ has already been defined, so one 
can project this direction out of the quantity which $B$ is applied to.
Thus, the algorithm can be modified to compute the vector
\begin{equation}
\label{eq:proj}
m_i \gets \theta_i \tilde u_{i-1} + B(w_{i-1} - \theta_i r_{i-1}), \quad
\theta_i \doteq \frac{\langle r,w \rangle}{||r||^2}
\end{equation}
Indeed, numerical experiments show that this modification can restore convergence in some cases.
This modification is, unfortunately, not useful as it stands for the purposes of improving the naive pipelined FCG algorithm, as it requires additional blocking reductions.
A natural question, however, is how one might estimate $\theta_i$.
One useful property is that as the preconditioner approaches $A^{-1}$, $\theta_i$ approaches $1$ for all $i$. 
This motivates the choice of $\theta_i \equiv 1$ as at least an ``unbiased'' approximation.

\subsubsection{Faithful Preconditioners}
To be more precise, we restrict our attention to a subclass of nonlinear preconditioners.
\begin{defn}
  A family of functions $B_n : \mathbb{C}^n \mapsto \mathbb{C}^n,\ n=1,2,\ldots$
  represents a \emph{faithful preconditioner} with respect to a family of linear 
  operators $A_n \in L(\mathbb{C}^n,\mathbb{C}^n)$ if $B_n(A_n(.)) - I_n(.)$ 
  is a  uniformly bounded sequence of operators, with respect to the standard norm on $\mathbb{C}^n$.
  That is, there exists $c \in \mathbb{R}$ such that
  \begin{equation*}
  ||B_n(A_n v) - v|| < c||v|| \quad \forall v \in \mathbb{C}^n \backslash \{0\},\quad
  \forall n \in \{1,2,\ldots\}
  \end{equation*}
  \revOne{19}{}
  \label{def:faithfulPC}
\end{defn}
\begin{rem} 
We will typically abuse notation and omit the dependence on $n$.\label{rem:notation}
\end{rem}
\begin{rem} 
Applied to the special case where $v = e_i \doteq x - x_i$,
 the error in a particular iteration of a Krylov method, noting that $B(Ae_i) = B(r_i) = u_i$,
  the condition above implies that 
\begin{equation*}
||u_i - e_i|| < c ||e_i||\text{,}
\end{equation*}
which is to say that the preconditioned residual is a \emph{uniformly accurate
approximation to the error}.
\label{rem:uminuse}
\end{rem}
\begin{rem}
\revTwo{6}{
The class of Sparse Approximate Inverse preconditioners are a special case of faithful preconditioners, if the approximation quality is uniform in the problem size.
}
\end{rem}

To compute $\theta_i$ in \eqref{eq:proj} requires additional inner products, so instead
we seek to cheaply estimate it. 
Note that $\theta_i$ minimizes the distance $||w_i - \theta_i r_i||$.
If we assume our preconditioner to be faithful with a small constant, so that the operator $AB(\cdot)-I(\cdot)$ is bounded with a small constant,
then $\theta_i\equiv 1$ 
is a reasonable estimate, as $||w - r||/||r|| = || A B(r) - r||/||r||$ is then small. 

Note that $\theta_i \equiv 0$ corresponds to the naive unrolling as in \eqref{eq:pcunroll2}.

\subsubsection{Defining an Improved Pipelined Flexible Conjugate Gradient Algorithm}

In a modified version of the naive pipelined Flexible Conjugate Gradients, 
Algorithm~\ref{alg:naivepipefcg},  we take advantage of the estimate $\theta_i \equiv 1$ to mitigate
the observed stagnation of convergence.

\label{sec:fcgmod}
Motivated by the considerations above, we observe great improvements employing the approximation
\begin{equation}
  \tilde u_i = B(r_i) = B(r_{i-1} - \alpha_{i-1}Ap_{i-1}) \approx  (1-\alpha_{i-1}) B(r_{i-1}) -\alpha_{i-1} B(Ap_{i-1} - r_{i-1})
  \label{eq:pcunroll3}
\end{equation}
instead of the recurrence relation~\eqref{eq:pcunroll2}.
This leads to a replacement of the algorithmic step \mbox{$m_i \gets B(w_i)$} in line
\ref{algline:naivepipefcgm} of Algorithm \ref{alg:naivepipefcg} by
\begin{equation}
  m_i \gets \tilde u_i + B(w_i - r_i).
  \label{eq:BwMinusR}
\end{equation}
For a linear preconditioner these are equivalent in exact arithmetic.

We proceed to show that this can be interpreted in terms of improved control of the accumulation of perturbations in $\tilde u_i$ with respect to $B(r_i)$,
assuming the preconditioner is sufficiently faithful. %Thought this is not strictly necessary - any linear part.

Note that $w_i = A \tilde u_i$ in exact arithmetic
and rewrite~\eqref{eq:BwMinusR} as
\begin{equation*}
m_i \gets \tilde u_i + B(A\tilde u_i - r_i) = \tilde u_i + B\left(A(\tilde u_i - e_i)\right).
\end{equation*}
With $q_i = m_i + \sum_k \beta_k q_k$ we write out the
recurrence relation
\begin{equation*}
\begin{split}
\tilde u_i = \tilde u_{i-1} - \alpha_{i-1} \left( \tilde u_{i-1} + B\big(A(\tilde u_{i-1}-e_{i-1}) \big) + \sum_k \beta_k q_k \right)\\
  = (1-\alpha_{i-1}) \tilde u_{i-1} - \alpha_{i-1} \left( B\big(A(\tilde u_{i-1}-e_{i-1}) \big) + \sum_k \beta_k q_k \right).
\end{split}
\end{equation*}
Rearranging terms we obtain
\begin{equation}
  \label{eq:uStabSplit}
  \tilde u_i  = \underbrace{\left[(1-\alpha_{i-1}) \tilde u_{i-1} - \alpha_{i-1}
    \sum_k \beta_k q_k\right]}_\text{I}\,
  - \,\underbrace{\left[\vphantom{\sum_k}\alpha_{i-1} B\big(A(\tilde u_{i-1}-e_{i-1}) \big)\right]}_\text{II},
\end{equation}
where \revTwo{7}{Term II} is affected by the nonlinear preconditioner but \revTwo{7}{Term I} is not.
Using Definition~\ref{def:faithfulPC} and Remark~\ref{rem:uminuse} we rewrite II as
\begin{equation*}
\alpha_{i-1} B\big(A(\tilde u_{i-1}-e_{i-1}) \big) = \alpha_{i-1} \tilde I \big(\tilde
u_{i-1}-e_{i-1}\big) \approx \alpha_{i-1} \big(\tilde u_{i-1}-e_{i-1}\big),
\end{equation*}
with the operator $\tilde I$ close to the identity operator using the notion of a faithful preconditioner. 
Comparing magnitudes of I and \revTwo{7}{this approximation of} II we expect that
\begin{equation}
  \label{eq:2}
  \|(1-\alpha_{i-1}) \tilde u_{i-1} - \alpha_{i-1}
    \sum_k \beta_k q_k\| > \| \alpha_{i-1} \big(\tilde u_{i-1}-e_{i-1}\big) \|
\end{equation}
as a faithful preconditioner transforms the residual such that it is similar to
the true error $e$. 
If instead we use $B(w)$ as in the naive Algorithm \ref{alg:naivepipefcg}, then
\begin{equation*}
\tilde u_i = \underbrace{\tilde u_{i-1} - \alpha_{i-1} \sum_k \beta_k q_k}_\text{III} - 
\underbrace{\vphantom{\sum_k} B\big(A\tilde u_{i-1}\big)}_\text{IV}.
\end{equation*}
Following the same line of arguments we obtain that $B(A\tilde u_{i-1})
= \tilde I (\tilde u_{i-1} ) \approx \tilde u_{i-1}$. Hence, Terms III and IV
are of similar magnitude.
While for the modified formulation the application of a nonlinear
preconditioner leads to small perturbations of the preconditioned residual and
consequently of the next search directions, this
does not hold true for the naive formulation. Through recurrence of the preconditioned
residual the accumulation of perturbations eventually leads to effectively random
search directions and stagnation of convergence.  
This is observed in Figure~\ref{fig:FCGtoy}, 
Once the ``noise level'' of the
preconditioner is reached, convergence stagnates.

By employing the modification to the computation of the quantity $m$, 
analyzed in the previous section, we arrive at
a useful pipelined variant of the FCG method, Algorithm \ref{alg:modpipefcg}.
%, also represented in Figure \ref{fig:modpipefcg}.
\begin{rem}
Algorithm \ref{alg:modpipefcg} satisfies several desiderata:
\ifx\undefined\nosiam\begin{romannum}\else\begin{itemize}\fi
\item In exact arithmetic, with a linear preconditioner, it is equivalent to CG and FCG.
\item As $B \to A^{-1}$, the operator becomes increasingly ``faithful'', hence the algorithm approaches the behavior of FCG.
\item No new reductions are introduced.
\item A single reduction phase per iteration is used, which can be overlapped with application of the operator and preconditioner.
\ifx\undefined\nosiam\end{romannum}\else\end{itemize}\fi
\end{rem}

The pipelined Flexible Conjugate Gradient method requires storing $(4\numax + 11)$ vectors, 
compared to $(\numax + 6)$ vectors for  FCG (Algorithm~\ref{alg:fcg}) and 
$(2\numax + 7)$ for Single Reduction FCG (Algorithm~\ref{alg:cgfcg}).

%%%%%%%%%%%%%%%%%%%%%%%%%%%%%%%%%%%%%%%%%%%%%%%%%%%%%%%%%%%%%%%%%%%%%%%%%%%%%%%
%%%%%%%%%%%%%%%%%%%%%%%%%%%%%%%%%%%%%%%%%%%%%%%%%%%%%%%%%%%%%%%%%%%%%%%%%%%%%%%
\section{Pipelined Generalized Conjugate Residuals Methods} 
%%%%%%%%%%%%%%%%%%%%%%%%%%%%%%%%%%%%%%%%%%%%%%%%%%%%%%%%%%%%%%%%%%%%%%%%%%%%%%%
%%%%%%%%%%%%%%%%%%%%%%%%%%%%%%%%%%%%%%%%%%%%%%%%%%%%%%%%%%%%%%%%%%%%%%%%%%%%%%%
\label{sec:cr}

The Conjugate Residual
family of Krylov subspace methods is applicable to systems with a Hermitian but 
not necessarily positive definite operator.
 While CG methods enforce $A$-conjugacy of search directions and minimize 
$||e||_A$, CR methods 
construct $A^HA$-orthogonal search directions.
In some derivations of the method, such as in Saad's textbook \cite{Saad2003}, 
iterates minimize $||e_i||_{A^HA} = ||r_i||$ over a Krylov subspace, which
puts these methods in close relation to MINRES \cite[\S6.4]{VanderVorst2003} and GMRES \cite{Saad1993,Saad1986}.
However, in some other derivations \cite{Luenberger:1970il}, these methods are described as
minimizing $||e||_A$, as does CG. %also see symmlq
Here, we consider methods of the first kind, that is minimal-residual methods.

%=============================================================================%
\subsection{Review of Conjugate Residuals Methods}
%=============================================================================%
\label{sec:review-conj-residuals}

There is also an ambiguity regarding the terminology
within the related flexible Krylov subspace methods. 
We note that the Generalized Conjugate
Residual (GCR) method \cite{Eisenstat:1983hk}, which this section is concerned with, is 
closely related to FCG.
However, despite its name, the GCR method is less closely related to Axelsson's
Generalized Conjugate Gradients (GCG)
method \cite{Axelsson:1991il,Axelsson:1987kc}. 
While FCG and GCR both
employ a multi term Gram-Schmidt orthogonalization of new search
directions against previous ones, GCG additionally applies old
directions in the update of the solution and residual vectors.

As in the preceding section, we describe the flexible method (GCR) 
and proceed to show its reduction to Preconditioned CR (PCR) for constant preconditioners
and extension to pipelined CR \cite{Ghysels2014}. 
We introduce pipelined GCR in
\S\ref{sec:pipelined-GCR}.

The GCR method is shown in Algorithm~\ref{alg:gcrR}. 
From a comparison of the GCR and FCG Algorithms~\ref{alg:gcrR} and \ref{alg:fcg}
it is straightforward to see that they employ $A^HA$- and $A$-inner products,
respectively.
As in Sec~\ref{sec:pcg},  we recover the non-flexible (or non-generalized)
Preconditioned Conjugate Residuals method by restricting to linear preconditioners.
In this case all $\beta_k$ but one are zero, giving Algorithm~\ref{alg:pcr}.

\ifx\undefined\nosiam\capstartfalse\fi % to get hypcap to ignore this wrapper figure*
\begin{figure*}
\begin{minipage}[t]{0.48\textwidth}
  \vspace{0pt}  
  \begin{algorithm}[H]
    \caption{Generalized CR}
    \label{alg:gcrR}
    \begin{algorithmic}[1]
      \FunctionDef{GCR}{$A$, $B$, $b$, $x_0$} 
      \StateNei{$r_0 \gets b - Ax_0$} 
      \StatePC{$u_0 \gets B(r_0)$}%this one doesn't actually need a tilde
      \StateLoc{$p_0 \gets u_0$}
      \StateNei{$s_0 \gets Ap_0$} 
      \StateRed{$\gamma_0 \gets \langle u_0,s_0 \rangle$}
      \StateRed{$\eta_0 \gets \langle s_0,s_0 \rangle$}
      \StateDef{$\alpha_0 \gets \gamma_0 / \eta_0$}
      \ForDef{$i = 1,2,\ldots$}
      \StateLoc{$x_i \gets x_{i-1} + \alpha_{i-1}p_{i-1}$}
      \StateLoc{$r_i \gets r_{i-1} - \alpha_{i-1}s_{i-1}$}\white{$()$}
      \StatePC{$u_i \gets B(r_i)$}
      \ForDef {$k=i-\nu_i,\ldots,i-1$}
      \StateRed{$\beta_{i,k} \gets \tfrac{-1}{\eta_k}\langle Au_i,s_k\rangle$}
      \EndFor
      \StateLoc{$p_i \gets u_i + \sum_{k=i-\nu_i}^{i-1}\beta_{i,k}p_k$}
      \StateNei{$s_i \gets Ap_i$}
      \StateRed{$\gamma_i \gets \langle u_i,s_i \rangle $}
      \StateRed{$\eta_i \gets \langle s_i,s_i\rangle$}
      \StateDef{$\alpha_i \gets \gamma_i / \eta_i$}
      \EndFor
      \EndFunction
    \end{algorithmic}
  \end{algorithm}
\end{minipage}%
\hspace{3px}
\begin{minipage}[t]{0.48\textwidth}
  \vspace{0pt}
  \begin{algorithm}[H]
    \caption{Preconditioned CR}
    \label{alg:pcr}
    \begin{algorithmic}[1]
      \FunctionDef{PCR}{$A$, $M^{-1}$, $b$, $x_0$} 
      \StateNei{$r_0 \gets b - Ax_0$}
      \StatePC{$u_0 \gets M^{-1}r_0$}
      \StateLoc{$p_0 \gets u_0$}
      \StateNei{$s_0 \gets Ap_0$}
      \StateRed{$\gamma_0 \gets \langle u_0,s_0 \rangle$}
      \StateRed{$\eta_0 \gets \langle s_0,s_0 \rangle$}
      \StateDef{$\alpha_0 \gets \gamma_0 / \eta_0$}
      \ForDef{$ i = 1,2,\ldots$}
      \StateLoc{$x_i \gets x_{i-1} + \alpha_{i-1}p_{i-1}$}
      \StateLoc{$r_i \gets r_{i-1} - \alpha_{i-1}s_{i-1}$}
      \StatePC{$u_i \gets M^{-1}r_i$}
      \StateRed{$\gamma_i \gets \langle u_i,s_i \rangle $}
      \StateDef{$\beta_i \gets \gamma_i / \gamma_{i-1}$}
        \vspace{4pt} % for alignment
      \StateLoc{$p_i \gets u_i + \beta_ip_{i-1}$}
      \StateNei{$s_i \gets Ap_i$}
      \StateDef{} % empty, for alignment
      \StateRed{$\eta_i \gets \langle s_i,s_i\rangle$}
      \StateDef{$\alpha_i \gets \gamma_i / \eta_i$}
      \EndFor
      \EndFunction
    \end{algorithmic}
  \end{algorithm}
\end{minipage}
\revTwo{9}{}
\end{figure*}
\ifx\undefined\nosiam\capstarttrue\fi % to get hypcap to ignore this wrapper figure*

%=============================================================================%
\subsection{Pipelined Generalized Conjugate Residual Methods}
%=============================================================================%
\label{sec:pipelined-GCR}

GCR and PCR do not allow any overlapping of global reductions in
their standard formulation. 
Reformulating the algorithms using the
Chronopoulos-Gear trick and introducing unrollings and pipelining
intermediates leads to two pipelined GCR variants shown in
algorithms~\ref{alg:PIPEGCR} and \ref{alg:PIPEGCR_w}.
These employ the same modification to the computation of $m_i$ as discussed in \S\ref{sec:fcgmod}.

In Algorithm~\ref{alg:PIPEGCR} the quantity $w_i$ is computed by performing the
sparse matrix-vector product $Au_i$ in line~\ref{algline:w_PIPEGCR}, whereas in
Algorithm~\ref{alg:PIPEGCR_w} $w$ is a recurred quantity. Note that the
overlapping properties of the variants are different as in
Algorithm~\ref{alg:PIPEGCR} the application of the preconditioner overlaps the
reductions while in Algorithm~\ref{alg:PIPEGCR_w} a sparse
matrix-vector product also contributes to the overlapping.
The performance of both variants is modeled in \S\ref{sec:perf-model}.
The suitability of one or the other algorithm
depends on the cluster architecture, in particular the speed of the connecting
network, and the available memory.
Algorithm~\ref{alg:PIPEGCR_w} has a larger memory footprint due to the
storage of the extra intermediate variable $z$ and its history of length $\nu_i$.

The original GCR method \cite{Eisenstat:1983hk} as presented in
Algorithm~\ref{alg:gcrR} requires a total of $(2\numax+7)$ vectors to be stored
%\footnote{Storing $\numax+1$ vectors $s$ can be avoided at the expense of $\numax+1$ sparse matrix-vector products $Ap_k$ in every iteration.})
. 
The pipelined variants of Algorithms~\ref{alg:PIPEGCR} and
\ref{alg:PIPEGCR_w} require $(\numax+2)$  and
$(2\numax+4)$  additional
vectors to be stored. \revTwo{10}{}
\ifx\undefined\nosiam\capstartfalse\fi % to get hypcap to ignore this wrapper figure*
\begin{figure*}
\begin{minipage}[t]{0.48\textwidth}
  \vspace{0pt}  
  \begin{algorithm}[H]
    \caption{Pipelined GCR \white{($w$ unrolled)}\label{alg:PIPEGCR}}
    \begin{algorithmic}[1]
      \FunctionDef{PIPEGCR}{$A$, $B$, $b$, $x_0$} 
      \StateNei{$r_0 \gets b - Ax_0$} 
      \StatePC{$\tilde u_0 \gets B(r_0)$}
      \StateLoc{$p_0 \gets u_0$}
      \StateNei{$s_0 \gets Ap_0$}
      \StatePC{$q_0 \gets B(s_0)$}
      \StateDef{\phantom.\white{$z_0 \gets Aq_0$}}
      \StateRed{$\gamma_0 \gets \langle r_0,Ar_0 \rangle$}
      \StateRed{$\eta_0 \gets ||s_0||$} \revOne{20}{}
      \StateDef{$\alpha_0 \gets \gamma_0 / \eta_0$}
      \ForDef{$ i = 1,2,\ldots$}
      \StateLoc{$x_i \gets x_{i-1} + \alpha_{i-1}p_{i-1}$}
      \StateLoc{$r_i \gets r_{i-1} - \alpha_{i-1}s_{i-1}$}
      \StateLoc{$u_i \gets \tilde u_{i-1} - \alpha_{i-1}q_{i-1}$}
      \StateNei{$w_i \gets Au_i$}\label{algline:w_PIPEGCR}
      \StatePC{$m_i \gets B(w_i-r_i) + \tilde u_i$}\label{algline:m_PIPEGCR}
      \StateDef{\phantom.\white{$n_i \gets Am_i$}}
      \StateRed{$\gamma_i \gets \langle r_i,w_i \rangle $}
      \StateRed{$\delta_i \gets \langle w_i,w_i \rangle$}
      \ForDef {$k = i-\nu_i,\ldots,i-1$}
      \StateRed{$\beta_{i,k} \gets \tfrac{-1}{\eta_k}\langle w_i,s_k\rangle$}
      \EndFor
      \StateLoc{$p_i \gets \tilde u_i + \sum_{k=i-\nu_i}^{i-1}\beta_{i,k}p_k$}
      \StateLoc{$s_i \gets w_i + \sum_{k=i-\nu_i}^{i-1}\beta_{i,k}s_k$}
      \StateLoc{$q_i \gets m_i + \sum_{k=i-\nu_i}^{i-1}\beta_{i,k}q_k$}
      \StateDef{\phantom.\white{$z_i \gets n_i + \sum_{k=i-\nu_i}^{i-1}\beta_{i,k}z_k$}}
      \StateDef{$\eta_i \gets \delta_i - \sum_{k=i-\nu_i}^{i-1}\beta^2_{i,k}\eta_k$}
      \StateDef{$\alpha_i \gets \gamma_i / \eta_i$}
      \EndFor
      \EndFunction
    \end{algorithmic}
  \end{algorithm}
\end{minipage}%
\hspace{3px}
\begin{minipage}[t]{0.5\textwidth}
  \vspace{0pt}
  \begin{algorithm}[H]
    \caption{Pipelined GCR ($w$ unrolled) \label{alg:PIPEGCR_w}}
    \begin{algorithmic}[1]
      \FunctionDef{PIPEGCR\_w}{$A$, $B$, $b$, $x_0$} 
      \StateNei{$r_0 \gets b - Ax_0$} 
      \StatePC{$\tilde u_0 \gets B(r_0)$}
      \StateLoc{$p_0 \gets \tilde u_0$}
      \StateNei{$ s_0 \gets Ap_0$}
      \StatePC{$q_0 \gets B(s_0)$}
      \StateNei{$z_0 \gets Aq_0$}
      \StateRed{$\gamma_0 \gets \langle r_0,Ar_0 \rangle$}
      \StateRed{$\eta_0 \gets ||s_0||$}
      \StateDef{$\alpha_0 \gets \gamma_0 / \eta_0$}
      \ForDef{$ i = 1,2,\ldots$}
      \StateLoc{$x_i \gets x_{i-1} + \alpha_{i-1}p_{i-1}$}
      \StateLoc{$r_i \gets r_{i-1} - \alpha_{i-1}s_{i-1}$}
      \StateLoc{$u_i \gets \tilde u_{i-1} - \alpha_{i-1}q_{i-1}$}
      \StateLoc{$w_i \gets w_{i-1} - \alpha_{i-1}z_{i-1}$}\label{algline:w_PIPEGCR_w}
      \StatePC{$m_i \gets B(w_i-r_i) + \tilde u_i$}
      \StateNei{$n_i \gets Am_i$}
      \StateRed{$\gamma_i \gets \langle r_i,w_i \rangle $}
      \StateRed{$\delta_i \gets \langle w_i,w_i \rangle$}
      \ForDef {$k = i-\nu_i,\ldots,i-1$}
      \StateRed{$\beta_{i,k} \gets \tfrac{-1}{\eta_k}\langle w_i,s_k\rangle$}
      \EndFor
      \StateLoc{$p_i \gets \tilde u_i + \sum_{k=i-\nu_i}^{i-1}\beta_{i,k}p_k$}
      \StateLoc{$s_i \gets w_i + \sum_{k=i-\nu_i}^{i-1}\beta_{i,k}s_k$}
      \StateLoc{$q_i \gets m_i + \sum_{k=i-\nu_i}^{i-1}\beta_{i,k}q_k$}
      \StateLoc{$z_i \gets n_i + \sum_{k=i-\nu_i}^{i-1}\beta_{i,k}z_k$}
      \StateDef{$\eta_i \gets \delta_i - \sum_{k=i-\nu_i}^{i-1}\beta^2_{i,k}\eta_k$}
      \StateDef{$\alpha_i \gets \gamma_i / \eta_i$}
      \EndFor
      \EndFunction
    \end{algorithmic}
  \end{algorithm}
\end{minipage}
\end{figure*}
\ifx\undefined\nosiam\capstarttrue\fi % to get hypcap to ignore this wrapper figure*

%%%%%%%%%%%%%%%%%%%%%%%%%%%%%%%%%%%%%%%%%%%%%%%%%%%%%%%%%%%%%%%%%%%%%%%%%%%%%
%%%%%%%%%%%%%%%%%%%%%%%%%%%%%%%%%%%%%%%%%%%%%%%%%%%%%%%%%%%%%%%%%%%%%%%%%%%%%
\section{Pipelined Flexible GMRES Methods}
%%%%%%%%%%%%%%%%%%%%%%%%%%%%%%%%%%%%%%%%%%%%%%%%%%%%%%%%%%%%%%%%%%%%%%%%%%%%%
%%%%%%%%%%%%%%%%%%%%%%%%%%%%%%%%%%%%%%%%%%%%%%%%%%%%%%%%%%%%%%%%%%%%%%%%%%%%%

\label{sec:gmres}
GMRES methods use explicitly-orthonormalized basis vectors and as such do not 
rely on the same identities which were so easily disturbed by algorithmic 
rearrangement in the FCG and GCR cases above.
Thus, deriving a usable pipelined FGMRES method
is a relatively straightforward extension of pipelined
GMRES, as predicted with the introduction of that method \cite{Ghysels2013}.

%=============================================================================%
\subsection{Notation}\label{sec:gmresnotation}
%=============================================================================%
We base our notation in Table \ref{tab:gmresnotation} on standard notation for
GMRES methods and pipelined GMRES.  
It is distinct from the notation used for CG and its variants
 because Flexible GMRES methods naturally use right preconditioning.
 The symbols $p_i$, $r_i$, $e_i$, $A$, $B$, and $M^{-1}$ have identical meaning to those in Table \ref{tab:cgnotation}.
\begin{table}[h!]
\footnotesize 
\begin{tabular}{l|l}
  $u_i$      & preconditioned basis vector $B(p_i)$ or $M^{-1}p_i$ \\
  $\tilde u$ & approximation to $u_i$, exact if $B$ is linear \\
  $z_i$ & transformed and preconditioned basis vector $Au_i$ or $A\tilde u_i$ \\
  $\bar z_i$ & preconditioned, transformed and shifted basis vector $Au_i -\sigma_ip_i$ or $A\tilde u_i-\sigma_ip_i$ \\ 
  $\bar q_i$ & pipelining intermediate $B(\bar z_i)$ \\
  $\bar w_i$ & pipelining intermediate $A\bar q_i$ \\
 \end{tabular}
\caption{Notation for GMRES and related methods}
\label{tab:gmresnotation}
\end{table}

%=============================================================================%
\subsection{Review of the Flexible GMRES Method}
%=============================================================================%
The Flexible GMRES\\ (\mbox{FGMRES}) method \cite{Saad1993} modifies the GMRES method \cite{Saad1986} to admit
variable right preconditioning. 
It accomplishes this by finding approximate solutions $x_i$ with $x_i-x_0$ having minimal residual 2-norm 
$||b-Ax_i||_2$ in the subspace $\text{span}(b,u_1,\ldots, u_{i-1})$,
 where $u_i \doteq B(p_i)$ and $p_i$ is formed by orthonormalizing
 $z_{i-1} \doteq  Au_{i-1} = AB(p_{i-1})$
with $p_0,\ldots,p_{i-1}$ \footnote{For improved numerical stability, the same subspace can be spanned by orthogonalizing 
$\bar z_{i-1} \doteq Au_{i-1} - \sigma_ip_{i-1}$ \cite{Ghysels2013}.}.
This Arnoldi process can be summarized as $AU = PH$, where $P$ is a matrix with $i$th column $p_i$,
$U$ is a matrix with $i$th column $u_i$, and $H$ is an upper Hessenberg matrix.
\revOne{21}{Let $P$ be a matrix with $i$th column $p_i$.} If $B = M^{-1}$ for some fixed linear operator,
GMRES is recovered: the solution is found in the Krylov space $\Kry^i(AM^{-1},b)$ and the relation
$AU=AM^{-1}P = PH$ is available to avoid the need to explicitly store the vectors
$u_i$. Algorithm \ref{alg:fgmres} describes the FGMRES method. Replacing
line~\ref{algline:fgmresupdatex} with
\begin{equation*}
x \gets x_0 + M^{-1}Py
\end{equation*}\revOne{21}{}
 in the case of a linear preconditioner recovers the right preconditioned GMRES algorithm.
The FGMRES method may be preferable even for a constant preconditioner if this
additional application of the preconditioner (``unwinding'') is expensive and the required additional memory is available.
Restarting and convergence checks (including on-the-fly QR decomposition of $H$ \revOne{22}{and happy breakdown checks}) are not included for simplicity. \revTwo{11}{Similarly, checks for the (rare in practice) possibility that the Hessenberg matrix $H$ is singular are omitted.}

\begin{algorithm}[!htbp]
\caption{Flexible GMRES \cite{Saad1993}}
\label{alg:fgmres}
\begin{algorithmic}[1]
\FunctionDef{FGMRES}{$A$,$B$,$b$,$x_0$,$m$}
\StateNei{$r_0 \gets b - Ax_0$}
\StateRed{$p_0 \gets r_0/||r_0||_2$}
\StatePC{$u_0 \gets B(p_0)$}
\StateNei{$z_0 \gets Au_0$}
\ForDef{$i=1,\ldots,m$}
\ForDef{$j=0,\ldots,i-1$}
\StateRed{$h_{j,i-1} \gets \langle p_j,z_{i-1}\rangle$}
\EndFor
\StateLoc{$\bar p \gets z_{i-1} - \sum_{k=0}^{i-1}p_j h_{j,i-1}$}
\StateRed{$h_{i,i-1} \gets ||\bar p||_2$}
\StateLoc{$p_i \gets \bar p / h_{i,i-1}$}
\StatePC{$u_i \gets B(p_i)$}
\StateNei{$z_i \gets Au_i$}
\EndFor
\StateLoc{$y \gets \text{argmin} ||H_{m+1,m}y - ||r_0||_2e_0||_2$}
\StateLoc{$x \gets x_0 + Uy$} \label{algline:fgmresupdatex}
\EndFunction
\end{algorithmic}
\end{algorithm}

%=============================================================================%
\subsection{Review of the Pipelined GMRES Method}
%=============================================================================%
The GMRES algorithm can be pipelined with the strategy outlined by Ghysels and Vanroose \cite{Ghysels2013}. 
They observe that at the cost of degraded numerical stability, 
GMRES can be modified to only involve one global reduction ($l^1$-GMRES in their nomenclature). 
This parallels the modification of preconditioned CG (Algorithm \ref{alg:pcg}) to derive Chronopoulos-Gear CG (Algorithm \ref{alg:cgcg}).
They use the identity $p_i = Ap_{i-1} - \sum_{k=1}^{i-1}\langle{Ap_{i-1},p_k}\rangle p_k$ 
%(where the $p_k$ are normalized) 
and note that if $z_i \doteq Ap_{i-1}$ is already available, this update can be performed while $Az_i = A^2p_{i-1}$ is being computed. 
Using the now available value of $p_i$, one can (locally) compute $z_{i+1} \doteq Ap_i$. 

%The pipelined GMRES method highlights the numerical difficulties inherent in the
%Gram-Schmidt process. It is well-known that the classical Gram-Schmidt method
%can be stabilized either by resorting to a modified Gram-Schmidt method, which
%does not parallelize well, or by employing iterative refinement,
%reorthogonalization, or even a restart of the algorithm. These difficulties
%increase with pipeline depth using the algorithm in \cite{Ghysels2013}.

%=============================================================================%
\subsection{Pipelining the Flexible GMRES Method} \label{sec:pipefgmres}
%=============================================================================%

The same procedure used to pipeline GMRES can also be used to produce a single-reduction variant of FGMRES, 
as shown in Algorithm~\ref{alg:cgfgmres}, and a pipelined version of FGMRES, shown in Algorithm~\ref{alg:naivepipefgmres}. 

This algorithm can suffer ``\revOne{2}{norm} breakdown'' requiring a
restart and refilling of the pipeline. 
Tuning of the shift parameters can help avoid this.
We allow for a distinct shift $\sigma_i$ to be used at each iteration,
but note that if this value is constant, the expression on line \ref{algline:pipefgmresz} 
of Algorithm \ref{alg:naivepipefgmres} simplifies to 
\begin{equation}\label{eq:pipefgmresz}
\bar z_i \gets (\bar w_{i-1} - \sum_{k=0}^{i-1}h_{k,i-1}\bar z_k)/h_{i,i-1}
\end{equation}
A practical choice, and the one we use in our implementations discussed in
\S~\ref{sec:implementation}-\ref{sec:tests}, is to set $\sigma_i$ to a constant value which approximates the 
largest eigenvalue of the preconditioned operator.

\begin{algorithm}[htbp]
\caption{Single Reduction FGMRES}
\label{alg:cgfgmres}
\begin{algorithmic}[1]
\FunctionDef{CGFGMRES}{$A$,$B$,$b$,$x_0$,$m$}
\StateNei{$r_0 \gets b - Ax_0$}
\StateRed{$p_0 \gets r_0/||r_0||_2$}
\StatePC{$u_0 \gets B(p_0)$}
\StateNei{$\bar z_0 \gets Au_0 - \sigma_0p_0$}
\ForDef{$i=1,\ldots,m$}
\ForDef{$j=0,\ldots,i-2$}
\StateRed{$\bar h_{j,i-1} \gets \langle p_j,\bar z_{i-1}\rangle$}
\StateDef{$h_{j,j-1} \gets \bar h_{j,i-1}$}
\EndFor
\StateRed{$\bar h_{i-1,i-1} \gets \langle p_{i-1},\bar z_{i-1}\rangle$}
\StateDef{$h_{i-1,i-1} \gets \bar h_{i-1,i-1} + \sigma_{i-1}$}
\StateRed{$t \gets ||\bar z_{i-1}||_2^2 - \sum_{k=0}^{i-1} \bar h_{k,i-1}^2$}
\IfDef{$t < 0$}
\StateDef{BREAKDOWN}
\EndIf
\StateDef{$\bar h_{i,i-1} \gets \sqrt{t}$}
\StateDef{$h_{i,i-1} \gets \bar h_{i,i-1}$}
\IfDef{$i = m$}
\Break
\EndIf
\StateLoc{$p_i \gets (\bar z_{i-1} - \sum_{k=0}^{i-1}\bar h_{k,i-1}p_k ) / \bar h_{i,i-1}$}
\StatePC{$u_i \gets B(p_i)$}
\StateNei{$\bar z_i \gets Au_i - \sigma_ip_i$}
\EndFor
\StateLoc{$y \gets \text{argmin} ||H_{m+1,m}y - ||r_0||_2e_0||_2$}
\StateLoc{$x \gets x_0 + Uy$} \label{algline:cgfgmresupdatex}
\EndFunction
\end{algorithmic}
\end{algorithm}

\begin{algorithm}[htbp]
\caption{Single-stage Pipelined FGMRES}
\label{alg:naivepipefgmres}
\begin{algorithmic}[1]
\FunctionDef{PIPEFGMRES}{$A$,$B$,$b$,$x_0$,$m$}
\StateNei{$r_0 \gets b - Ax_0$}
\StateRed{$p_0 \gets r_0/||r_0||_2$}
\StatePC{$u_0 \gets B(p_0)$}
\StateNei{$\bar z_0 \gets Au_0 - \sigma_0p_0$}
\StatePC{$\bar q_0 \gets B(\bar z_0)$}
\StateNei{$\bar w_0 \gets Aq_0$}
\ForDef{$i=1,\ldots,m$}
\ForDef{$j=0,\ldots,i-2$}
\StateRed{$\bar h_{j,i-1} \gets \langle p_j,\bar z_{i-1}\rangle$}
\StateDef{$h_{j,j-1} \gets \bar h_{j,i-1}$}
\EndFor
\StateRed{$\bar h_{i-1,i-1} \gets \langle p_{i-1},\bar z_{i-1}\rangle$}
\StateDef{$h_{i-1,i-1} \gets \bar h_{i-1,i-1} + \sigma_{i-1}$}
\StateRed{$t \gets ||\bar z_{i-1}||_2^2 - \sum_{k=0}^{i-1} \bar h_{k,i-1}^2$}
\IfDef{$t < 0$}
\StateDef{BREAKDOWN}
\EndIf
\StateDef{$\bar h_{i,i-1} \gets \sqrt{t}$}
\StateDef{$h_{i,i-1} \gets \bar h_{i,i-1}$}
\IfDef{$i = m$}
\Break
\EndIf
\StateLoc{$p_i \gets (\bar z_{i-1} - \sum_{k=0}^{i-1}\bar h_{k,i-1}p_k ) / \bar h_{i,i-1}$}
\StateLoc{$\tilde u_i \gets (\bar q_{i-1} - \sum_{k=0}^{i-1}u_k \bar h_{k,i-1})/\bar h_{i,i-1}$} \label{algline:naivepipefgmresupdateu}
\StateLoc{$\bar z_i \gets (\bar w_{i-1} - \sum_{k=0}^{i-1}(\bar z_k + \sigma_kp_k) \bar h_{k,i-1})/\bar h_{i,i-1} - \sigma_i p_i$} 
\Comment see \eqref{eq:pipefgmresz} \label{algline:pipefgmresz}
\StatePC{$\bar q_i \gets B(\bar z_i)$}
\StateNei{$\bar w_i \gets A\bar q_i$}
\EndFor
\StateLoc{$y \gets \text{argmin} ||H_{m+1,m}y - ||r_0||_2e_0||_2$}
\StateLoc{$x \gets x_0 + \tilde Uy$} \label{algline:naivepipefgmresupdatex}%Again, notation not totally clear
\EndFunction
\end{algorithmic}
\end{algorithm}

In exact arithmetic, the iterates produced by Algorithm~\ref{alg:naivepipefgmres} 
are not equivalent to the FGMRES algorithm because of the approximation in
computing $\tilde u_i$ on line \ref{algline:naivepipefgmresupdateu}. 
Recall that $u_i \doteq B(p_i) = B((v_i - \sum_{k=1}^{-1}h_{k,i-1}p_k)/h_{i,i-1})$. 
If $B$ were linear, then the approximation $\tilde u_i$ on line \ref{algline:naivepipefgmresupdateu} 
would be arithmetically equivalent.
PIPEFGMRES requires storage of $4m + 2$ vectors, as opposed to $2m + 2$ for FGMRES and CGFGMRES.

%%%%%%%%%%%%%%%%%%%%%%%%%%%%%%%%%%%%%%%%%%%%%%%%%%%%%%%%%%%%%%%%%%%%%%%%%%%%%%%
%%%%%%%%%%%%%%%%%%%%%%%%%%%%%%%%%%%%%%%%%%%%%%%%%%%%%%%%%%%%%%%%%%%%%%%%%%%%%%%
\section{Implementation} \label{sec:implementation}
%%%%%%%%%%%%%%%%%%%%%%%%%%%%%%%%%%%%%%%%%%%%%%%%%%%%%%%%%%%%%%%%%%%%%%%%%%%%%%%
%%%%%%%%%%%%%%%%%%%%%%%%%%%%%%%%%%%%%%%%%%%%%%%%%%%%%%%%%%%%%%%%%%%%%%%%%%%%%%%

As a practical consequence of our choice of pipelined algorithms which are not 
arithmetically equivalent to their base flexible methods, useful preconditioners
for the base solver, for a given problem, will not necessarily translate to 
be effective solvers with the pipelined variant.
In particular, for the Pipelined FCG and GCR methods, the preconditioner must be compatible with
the approximation $\theta_i \equiv 1$, as described in \S\ref{sec:fcg-issues}.
Thus, experimentation is required, as it is in general when choosing preconditioners for
complex systems where canonical solutions or analysis tools are not available.

To allow for this experimentation to be performed easily, and to allow for
usage in highly distributed settings by way \revOne{23}{of} MPI, we implement 
the PIPEFCG (Alg.~\ref{alg:modpipefcg}),
PIPEGCR (Alg.~\ref{alg:PIPEGCR}),
 PIPEGCR\_w (Alg.~\ref{alg:PIPEGCR_w}), and PIPEFGMRES (Alg.~\ref{alg:naivepipefgmres}) algorithms
within the \texttt{KSP} class in \textsc{PETSc} \cite{petsc-user-ref,petsc-web-page}. 
The implementations are identified as \texttt{KSPPIPEFCG}, \texttt{KSPPIPEGCR}, and \texttt{KSPPIPEFGMRES}\revTwo{17}{}\footnote{See \texttt{bitbucket.org/pascgeopc/petsc} for current information.}.

An important practical point is that although many residual norms may be
used to monitor convergence, only the natural norm for the CG- and CR-based solvers
may be used without introducing extra computations and reductions.
Similarly, the right preconditioned residual norm for PIPEFGMRES
should be used when performance is important. 
Other norms can and should be used to assess convergence behavior.

All of the pipelined algorithms may suffer from norm breakdown, although in practice this is typically
only observed near convergence. 
Our implementations detect norm breakdown and respond by restarting the algorithm, flushing the pipeline in the process.

\begin{rem}[Truncation strategies for FCG and GCR methods]
In \cite{Notay2000} Notay proposes a truncation--restart technique for determining 
the number of old directions $\nu_i$ to be used in the current iteration 
$i$ as $\nu_i = \max \left( 1,\mod(i,\numax+1) \right)$. 
With this technique after every restart there are two consecutive iterations
using only the respective last search direction. Consequently, the first two
steps after a restart are actually standard CG iterations. 
In our implementations we implement the rule $\nu_i = \mod(i,\numax) + 1$. This
employs $n+1$ directions in the $n$th iteration after a restart.
If $\numax = 1$ we recover the IPCG algorithm as analyzed by Golub and Ye \cite{Golub1999}.
We also provide a standard truncation procedure, $\nu_i = \min(i,\numax)$.
\end{rem}

\begin{rem}[Shift parameters for \texttt{KSPPIPEFGMRES}]
Our implementation only covers the constant-$\sigma$ case, as discussed in \S\ref{sec:pipefgmres}.
\end{rem}

\begin{rem}[Scaling] 
We observe that scaling the system matrix $A$ is
often essential to allow for a faithful preconditioner (and see the survey by Wathen \cite{Wathen2015}
for an exposition on scenarios which may otherwise profit from diagonal scaling). 
This is due to the fact that an effective preconditioner for $A$ may not be faithful, 
as we have defined it in Definition \ref{def:faithfulPC}.
Indeed, a Krylov method will typically converge rapidly when the preconditioned operator
has a small number of tight clusters of eigenvalues, a property which is not affected by scaling of the preconditioner or operator. 
Operators close to the identity frequently exhibit this property, but this property is neither necessary nor sufficient.
%A scaling of the operator may be able to make an effective preconditioner faithful.
\end{rem}
%%%%%%%%%%%%%%%%%%%%%%%%%%%%%%%%%%%%%%%%%%%%%%%%%%%%%%%%%%%%%%%%%%%%%%%%%%%%%%%
%%%%%%%%%%%%%%%%%%%%%%%%%%%%%%%%%%%%%%%%%%%%%%%%%%%%%%%%%%%%%%%%%%%%%%%%%%%%%%%
\section{Numerical Experiments} \label{sec:tests}
%%%%%%%%%%%%%%%%%%%%%%%%%%%%%%%%%%%%%%%%%%%%%%%%%%%%%%%%%%%%%%%%%%%%%%%%%%%%%%%
%%%%%%%%%%%%%%%%%%%%%%%%%%%%%%%%%%%%%%%%%%%%%%%%%%%%%%%%%%%%%%%%%%%%%%%%%%%%%%%

% PDS: could put something here if you don't like having two headings one after
%  the other

%=============================================================================%
\subsection{Blocking and Non-Blocking Reductions}
%=============================================================================%
\label{sec:comp-block-non-block}

%PDS: I ran some new daint tests with the permutron (in the data scratch repo and processed with a google sheet). 
%     Those can be worked in here as well to give some idea of the variability on 1024 nodes.

The MPI-3 standard \cite{MPI3} provides an API for asynchronous operations which
\emph{may} allow for overlap of communication and computation, as required by
the algorithms described in this paper. 
This is still an active area of software development \cite{Bruggencate2010,Hoefler2008b,Wittmann2013}.
The latency involved in a global reduction on a given massively parallel 
machine is difficult to characterize, particularly if only using a portion of the nodes, 
and performance models are not always predictive \cite{Hoefler2008}. 
It should be noted that the only reliably reproducible quantity when testing message 
passing routines on complex machines is the \emph{minimum} time required for an operation, 
in the limit of a large number of tests,
though of course users are more affected by typical or average times \cite{Gropp1999}.

Here and in the following tests, we report a variety of data and statistics to allow for \emph{interpretability} \cite{Hoefler2015} of our experiments,
even in the absence of access to the specific hardware and software environment available to
the authors at the time of this writing. Note that the performance model discussed in \S\ref{sec:perf-model} can
also provide additional insight into the performance characteristics of the new solvers.

In order to assess the potential time savings and performance enhancements by using 
pipelined versus non-pipelined preconditioned Krylov solvers, 
empirical data are presented, comparing the time to execute
\texttt{MPI\_Allreduce} and \texttt{MPI\_Iallreduce} \cite{MPI3}. 
The reductions are paired with identical computational work local to each rank, \revOne{24}{consisting of arbitrary floating point operations on individual elements of a small local array}. 
In the former case we get the combined time of the
sequence of local work and the blocking reduction while in the latter case, we
time the sequence of starting the reduction, doing the local work and finishing
with \texttt{MPI\_Wait}. The reduction size is 32 double
values. Table~\ref{tab:reduction_timing} lists the timing results for
Piz Daint, a Cray XC30 cluster at the Swiss National Supercomputing Center (CSCS).
The machine includes 5272 nodes, 
each with an 8-core Intel Xeon E5-2670 CPU, 32 GiB DDR3-1600 RAM, and
Aries routing and communications ASIC.
The cluster features the Dragonfly network topology \cite{Kim:eg}.
%\texttt{hexagon}: Cray XE6, 696 nodes, 2×16 AMD Interlagos CPUs, Gemini 2.5D
%torus network hosted by the Bergen Centre for Computational Science (BCCS), and
%supported by Notur -- the Norwegian Metacenter for Computational Science
As a reference we also include
timings for performing the reduction and the local work only. We used five warmup rounds and
50 trial runs. The timings presented are the minimum and maximum times across
ranks averaged over the trial runs.

It can be confirmed that \texttt{MPI\_Iallreduce}
functionality is in place.
While the time for performing the
local work and a blocking reduction just exceeds the cumulative time for the
individual operations, the reduction time is hidden when using non-blocking
reductions. From the min/max values we also see that usually timing is
consistent across the communicator sizes we considered.

Table~\ref{tab:reduction_timing} is also instructive to understand when a
performance gain can be expected by using a pipelined solver as these
methods require computing and updating of additional intermediate variables.

\begin{table}[h!]
  \begin{center}
    \begin{tabular}{l|rrrrrrrrr}
      \toprule
      \#Ranks & 256 & 1024 & 2048 & 4096 & 8192 & 16384 & 24576 & 32768 \\
      \midrule
      \multirow{2}{*}{Reduction only}  & 0.188 & 0.577  & 0.613 & 0.670 &  0.787  & 0.843 & 0.918  & 0.894 \\
              & 0.198 & 0.588 & 0.622 & 0.689 &  0.821 & 0.876& 0.956 & 0.929\\
      \midrule
      \multirow{2}{*}{Local work only} & 4.88 &4.88 &4.88 &4.87 &4.88 &4.88 &4.88 &4.88 \\
              & 4.95&4.95&4.97&4.97&4.98&5.68&5.67&5.68\\
      \midrule
      \multirow{2}{*}{+ blocking red.} & 5.56&5.75 &5.81 &5.89 &6.03 &6.12 &6.13 &6.13\\
              & 5.57&5.76&5.83&5.90&6.05&6.14&6.15&6.15\\
      \midrule
      \multirow{2}{*}{+ non-bl. red.} & 4.90 &4.88 &4.90&4.89&4.93&5.65&5.64&5.66 \\
              & 4.98&4.98&4.99&4.99&4.98&5.71&5.69&5.71\\
      \bottomrule
      %\multicolumn{9}{c}{hexagon}\\
      %\hline
      %\multirow{2}{*}{Reduction only} &\multirow{2}{*}{--}&4.23 &4.23 &4.63 &\multirow{2}{*}{--}&\multirow{2}{*}{--}&\multirow{2}{*}{--}&\multirow{2}{*}{--}\\
      %        &&4.49&4.50&4.92&&&&\\
      %\hline
      %\multirow{2}{*}{Local work only (LW)} &\multirow{2}{*}{--}&13.7&13.6 &13.6 &\multirow{2}{*}{--}&\multirow{2}{*}{--}&\multirow{2}{*}{--}&\multirow{2}{*}{--}\\
      %        &&15.1&15.0&15.1&&&&\\
      %\hline
      %\multirow{2}{*}{LW + blocking red.} &\multirow{2}{*}{--}&20.2 &20.6 &20.7 &\multirow{2}{*}{--}&\multirow{2}{*}{--}&\multirow{2}{*}{--}&\multirow{2}{*}{--}\\
      %        &&20.5&20.9&21.0&&&&\\
      %\hline
      %\multirow{2}{*}{LW + non-bl. red.} &\multirow{2}{*}{--}&14.3&14.3&14.4&\multirow{2}{*}{--}&\multirow{2}{*}{--}&\multirow{2}{*}{--}&\multirow{2}{*}{--}\\
      %        &&14.5&14.5&14.6&&&&\\
      %\hline\hline
    \end{tabular}
    \caption{Timing for blocking and non-blocking reductions in milliseconds. We report both the
      min (upper rows) and max (lower rows) CPU time over all ranks within the MPI communicator. 
      The times for the non-blocking reductions with local
      work have to be compared to the maximum times of the local work only as
      the reduction can terminate only once all ranks finish with the local work.}
    \label{tab:reduction_timing}
  \end{center}
\end{table}

%=============================================================================%
\subsection{3D Strong-Scaled Coarse Grid Solver Example}
%=============================================================================%

A key motivation for the development of pipelined Krylov methods for use on 
current supercomputers is the scenario in which the amount of local work per 
iteration in a Krylov method is very low, hence requiring time comparable to that consumed by reductions. 
This is the case when attempting to strong-scale past the point of the current reduction bottleneck. 
The same situation arises when a Krylov method is employed as a coarse-grid solver within a multigrid hierarchy.

In the Piz Daint environment described in \S\ref{sec:comp-block-non-block}, we use
Cray CLE 5.2.40 and Cray-MPICH 7.2.2,
a \textsc{PETSc} development branch based on \textsc{PETSc} 3.6 \texttt{maint}\footnote{commit \texttt{2c5c660747c89d9a265191735ba48020b7a0dd72} at \texttt{https://bitbucket.org/petsc/petsc}},
and a \textsc{pTatin3D} \cite{May2014} development branch. %updated for the \textsc{PETSc} 3.6 API.
This allows us to evaluate our solvers with real application software to solve 
the variable viscosity Stokes equations using $\mathbb{Q}_2-\mathbb{P}_1^{\text{disc}}$ mixed finite elements.

We examine a 3D ``viscous sinker'' scenario with a spherical inclusion 
of higher ($1.2\times$) density and viscosity ($100\times$) in a cubic domain.
This is a useful test case for challenging Stokes flow problems, as it presents a 
tunable, non-grid-aligned viscosity contrast.
Tests correspond to a highly-distributed coarse grid solve, or an extreme strong-scaling.
The former case is a particularly interesting potential application of our methods, 
as it allows some adjustment of the grid size to balance communication and computation to be overlapped.
$4096$ MPI ranks on $1024$ nodes each hold a single
$\mathbb{Q}_2$ finite element, for $16^3$ elements and $107,811$ degrees of freedom.
The resulting linear systems are solved with the new \texttt{KSPFCG} and \texttt{KSPPIPEFCG} solvers to a relative tolerance of $10^{-8}$ in the natural residual norm.
In each case, $\nu_\text{max}=30$ and standard truncation is employed.
FCG and PIPEFCG are compared for two nonlinear preconditioners.
60 tests were run for each of the two solves in each test, 10 at a time per batch job. 
The first run is considered warmup and is reported separately.
Sets of three batch jobs were submitted at a time.
\subsubsection{Block Jacobi Preconditioner}
\begin{figure}[!htbp]
\begin{center}
\includegraphics[width=0.49\textwidth]{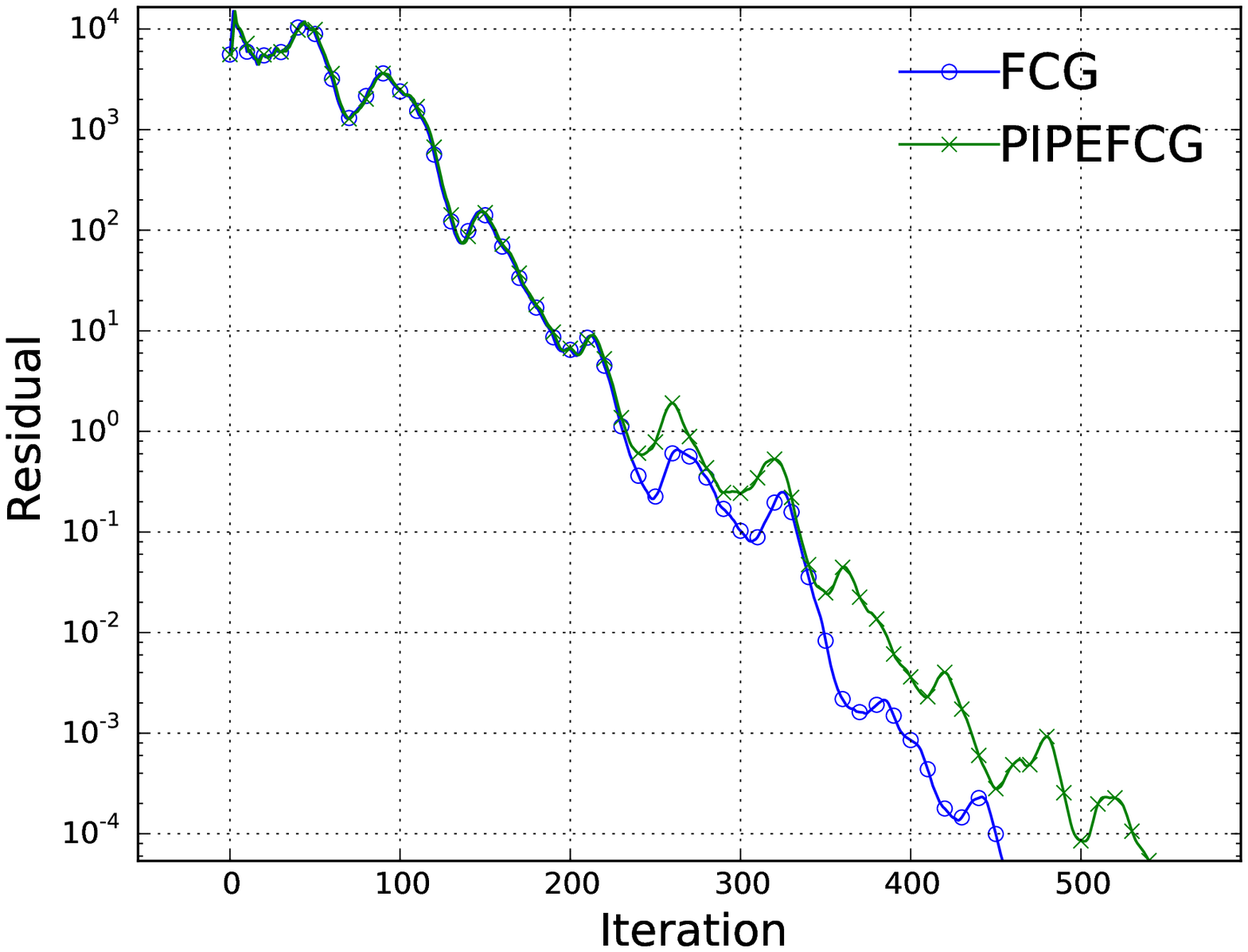}\hfill
\includegraphics[width=0.49\textwidth]{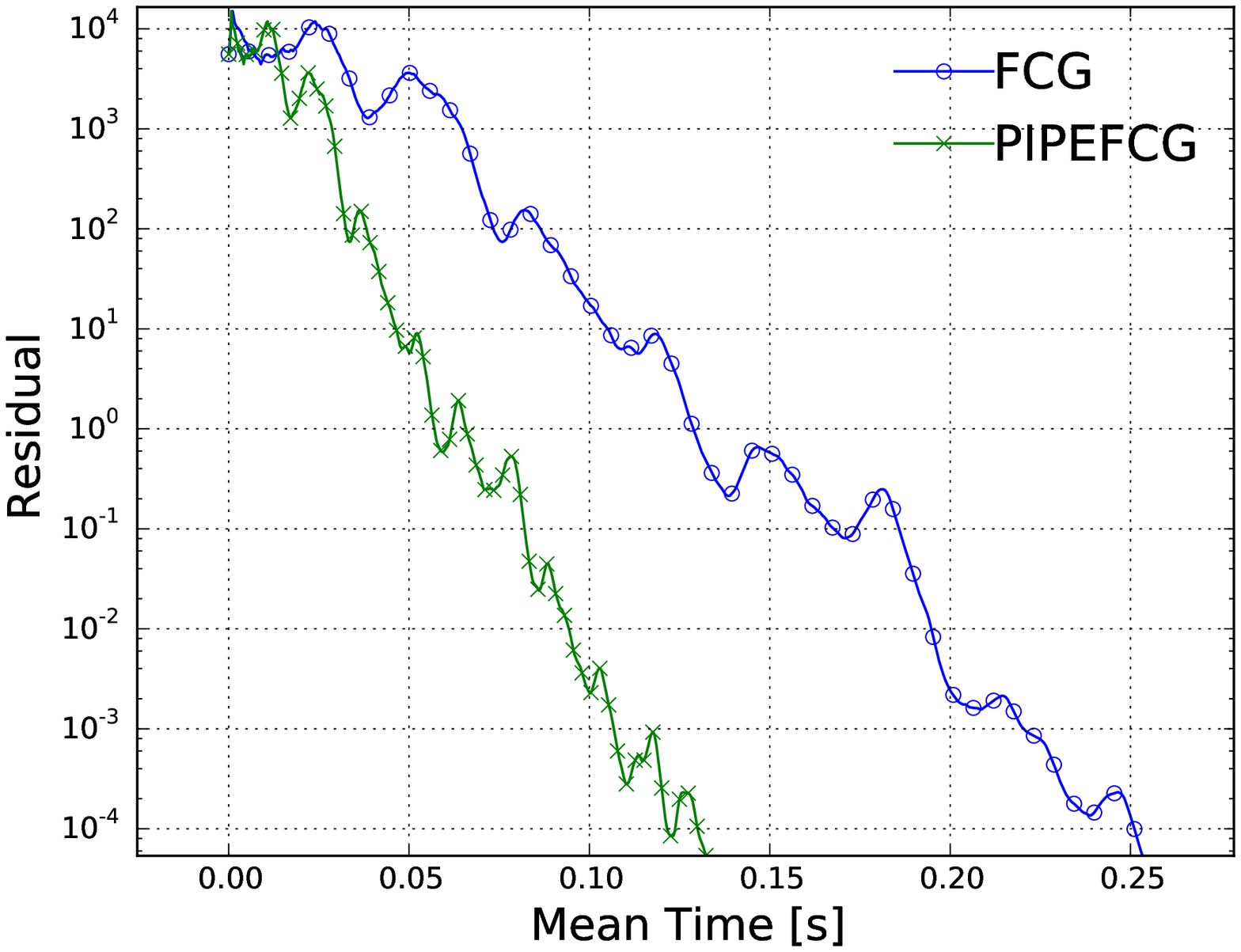}
% Table goes here
\caption{Convergence behavior for the viscous sinker system using block Jacobi preconditioning\label{fig:bjac}.}
\end{center}
\end{figure}

A per-rank block Jacobi preconditioner is employed, performing 5 iterations of Jacobi-preconditioned CG on each block.
Figure \ref{fig:bjac} shows the residual norm convergence behavior as a function of the iteration count (left) and CPU time (right) for a typical experiment. 
Despite the increase in iterations required for convergence when using PIPEFCG, the time to solution is $\sim 2\times$ faster compared to FCG.
Table \ref{tab:bjacresults} collects statistical variations obtained from multiple runs 
(indicated via the ``Samples'' column). 
The ``First solve'' column indicates whether the solve was the first call to 
\textsc{PETSc}'s \texttt{KSPSolve} routine after defining the system; 
successive solves were run with no intermediate calls apart from those used for timing.
Note the relatively large variability in timings. 
Interestingly, the maximum speedup per iteration can \emph{exceed} $2\times$. 
While this might initially seem impossible to accomplish with an algorithm which aims to 
overlap two processes, it should be noted that there is also some benefit to be gained in 
relaxing the synchronization required, beneficial in the presence of heavy-tailed 
distributions of local processing times. 
For more, see the recent work by Knepley et al. \cite{Morgan2015}.
\begin{table}[htbp]
\begin{center}
\footnotesize
%  \begin{tabular}{lcccccccccc}
%  \hline 
%         &             &         & \multicolumn{4}{c}{Total Linear Solve Time [s]}       \\
%  Solver & First solve & Samples & Mean (std. dev.)        & Min      & Median   & Max   \\
%  \hline
%  FCG     & yes        & 6       & 2.90E-01 (6.27E-02) & 2.04E-01 & 2.92E-01 & 3.58E-01  \\
%  PIPEFCG & yes        & 6       & 2.27E-01 (3.89E-02) & 1.79E-01 & 2.23E-01 & 2.92E-01  \\
%  FCG     & no         & 54      & 2.53E-01 (5.03E-02) & 1.70E-01 & 2.49E-01 & 4.15E-01  \\
%  PIPEFCG & no         & 54      & 1.32E-01 (1.70E-02) & 1.05E-01 & 1.30E-01 & 2.01E-01  \\
%  \hline
%  \end{tabular}
%  \vspace{10px}
  \begin{tabular}{lll l lll}
  \toprule 
         &             &         & &  \multicolumn{3}{c}{Time / Krylov Iteration [s]}\\
  Solver & First solve & Samples &Its. & Mean (Std. Dev.)        & Min.     & Max.   \\
  \midrule
  FCG         & yes        & 6   &454    & 6.38E-04 (1.38E-04) & 4.48E-04 & 7.89E-04  \\
  PIPEFCG & yes        & 6   &540    & 4.21E-04 (7.20E-05) & 3.31E-04 & 5.41E-04  \\
  \midrule
  FCG         & no         & 54   &454   & 5.58E-04 (1.11E-04) & 3.73E-04 & 9.13E-04  \\
  PIPEFCG & no         & 54   &540   & 2.45E-04 (3.15E-05) & 1.94E-04 & 3.73E-04  \\
  \bottomrule
  \end{tabular}
  \caption{Statistics for multiple runs of the \textsc{pTatin3D} viscous sinker system with a block Jacobi preconditioner. 
   \label{tab:bjacresults}}
\end{center}
\end{table}

\subsubsection{RASM-1 Preconditioner}
\begin{figure}[!htbp]
\begin{center}
\includegraphics[width=0.49\textwidth]{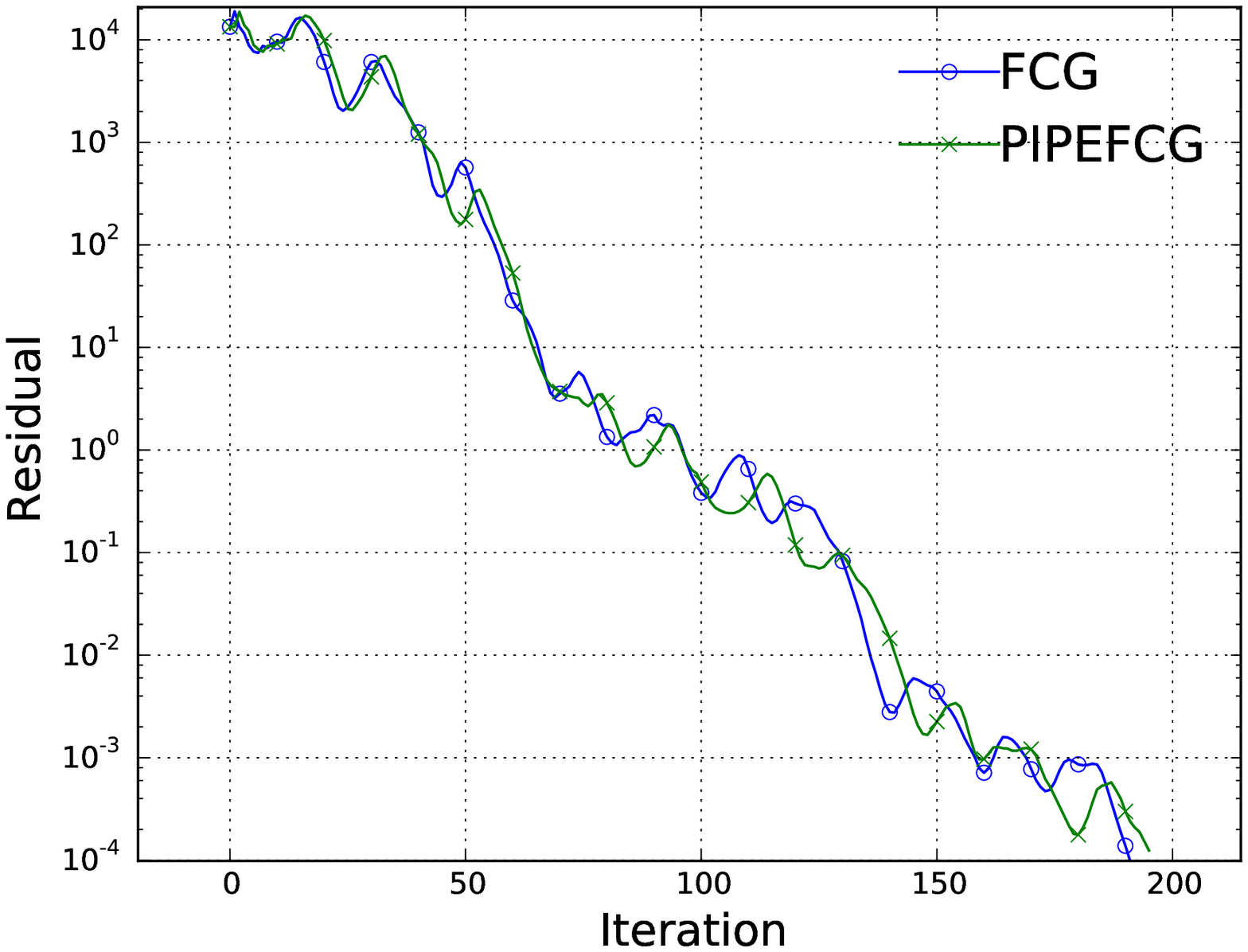}\hfill
\includegraphics[width=0.49\textwidth]{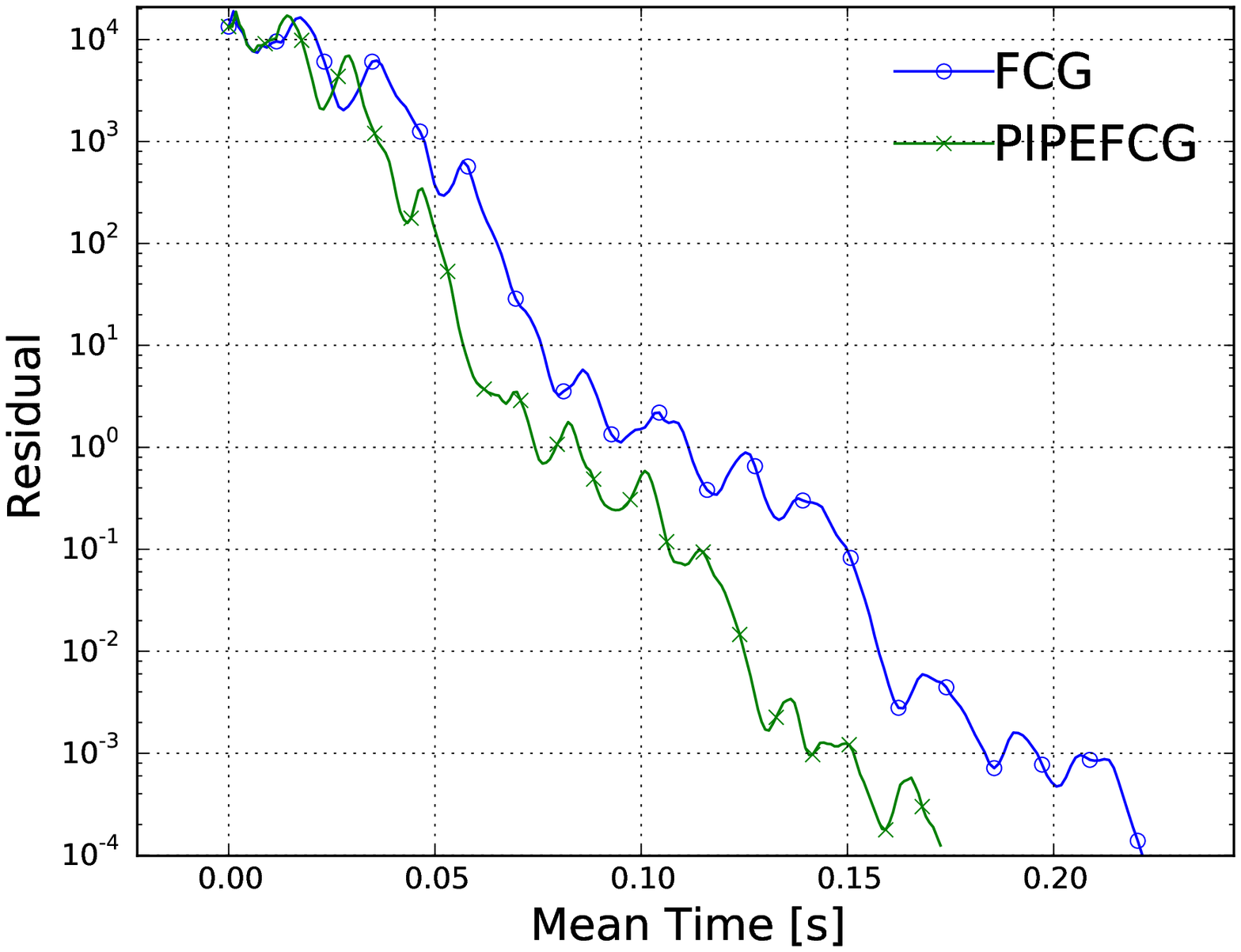}
% Table goes here
\caption{Convergence behavior for the viscous sinker system using RASM-1 preconditioning\label{fig:asm}}
\end{center}
\end{figure}

A rank-wise RASM-1 \cite{Cai1999} preconditioner is used, with 5 iterations of Jacobi-preconditioned CG as a subdomain solver.
Figure \ref{fig:asm} shows the residual norm convergence behavior as a function of the iteration count (left) and mean CPU time (right) for a typical experiment. 
Note that, as expected, choosing a more faithful preconditioner allows recovery of comparable convergence behavior between the pipelined and base method. 
Here, the speedup is not as pronounced as in the previous example, as the problem size was 
tuned as an example of a coarse grid solve with respect to the block Jacobi preconditioner.
This example shows that speedup is practically achievable even given nested \textsc{MPI} calls defining neighbor-wise communication patterns 
in the preconditioner.
Table \ref{tab:asmresults} presents statistics of the multiple solves performed. 
\begin{table}[H]
\begin{center}
\footnotesize
%  \begin{tabular}{lcccccc}
%  \hline 
%         &             &         & \multicolumn{4}{c}{Total Linear Solve Time [s]}       \\
%  Solver & First solve & Samples & Mean (std. dev.)        & Min      & Median   & Max       \\
%  \hline
%  FCG     & yes        & 6       & 2.44E-01 (2.45E-02) & 2.12E-01 & 2.48E-01 & 2.70E-01  \\
%  PIPEFCG & yes        & 6       & 2.36E-01 (2.29E-02) & 2.14E-01 & 2.28E-01 & 2.64E-01  \\
%  FCG     & no         & 54      & 2.22E-01 (2.83E-02) & 1.73E-01 & 2.15E-01 & 2.97E-01  \\
%  PIPEFCG & no         & 54      & 1.73E-01 (1.93E-02) & 1.36E-01 & 1.73E-01 & 2.15E-01  \\
%  \hline
%  \end{tabular}
%  \vspace{10px}
  \begin{tabular}{llll lll}
  \toprule 
         &             &         & &  \multicolumn{3}{c}{Time / Krylov Iteration [s]}     \\
  Solver & First solve & Samples & Its. & Mean (Std. Dev.)        & Min.      & Max.       \\
  \midrule
  FCG          & yes        & 6   &191    & 1.28E-03 (1.28E-04) & 1.11E-03 & 1.41E-03  \\
  PIPEFCG  & yes        & 6   &195    & 1.21E-03 (1.18E-04) & 1.10E-03 & 1.36E-03  \\
  \midrule
  FCG          & no         & 54  &191    & 1.16E-03 (1.48E-04) & 9.08E-04 & 1.55E-03   \\
  PIPEFCG  & no         & 54  &195    & 8.85E-04 (9.89E-05) & 6.99E-04 & 1.10E-03  \\
  \bottomrule 
  \end{tabular}
  \caption{Statistics for multiple runs of the \textsc{pTatin3D} viscous sinker system with restricted ASM-1 preconditioner. 
  \label{tab:asmresults}}
\end{center}
\end{table}

%%=============================================================================%
\subsection{2D Stokes Flow Test}
%%=============================================================================%
\label{sec:ex43}
As a relevant yet still comparatively simple test case of a full solve, we ran extensive tests
using a \textsc{PETSc} tutorial\footnote{\texttt{www.mcs.anl.gov/petsc/petsc-current/src/ksp/ksp/examples/tutorials/ex43.c.html}},
which solves the incompressible, variable viscosity Stokes equation in two
spatial dimensions. 
The discretization employs $\mathbb{Q}_1-\mathbb{Q}_1$ elements, stabilized with
Bochev's polynomial projection method \cite{Dohrmann:2004er} and
free slip boundary conditions on all sides. The viscosity structure is a circular viscous sinker,
with density contrast giving a buoyant force and a viscosity contrast of $25\times$.
All tests run on 4096 MPI ranks on 1024 nodes of Piz Daint.
We experiment with all three new solvers as subsolvers within an upper-triangular block preconditioner.
%The preconditioned system is of the form
%\begin{equation*}
%\left[ \begin{array}{cc} K & G \\ D & C \end{array}\right] \left[ \begin{array}{cc}\hat K & G \\ 0 & \hat S \end{array}\right]^{-1} = \left[\begin{array}{c} f \\ 0 \end{array}\right] 
%\end{equation*}
The discrete Stokes problem is of the form 
\begin{equation*}
\left[ \begin{array}{cc} K & G \\ D & C \end{array}\right] \left[ \begin{array}{c}u \\ p \end{array}\right] = \left[\begin{array}{c} f \\ 0 \end{array}\right] 
\end{equation*}
which is right preconditioned with
\begin{equation*}
B = \left[ \begin{array}{cc}\hat K & G \\ 0 & \hat S \end{array}\right]^{-1},  
\end{equation*}
where $\hat S$ is a Jacobi preconditioner derived from the pressure mass matrix scaled by the local (element-wise) inverse of the viscosity
and $\hat K$ is an inexact flexible pipelined Krylov solve applied to the viscous block of the Stokes system, preconditioned with a rank-wise 
block Jacobi preconditioner with 5 iterations of Jacobi-preconditioned CG as subsolves.
The convergence criterion is that the unpreconditioned norm of an outer FGMRES method is reduced by a factor of $10^6$.
Inner solves are limited to 100 iterations.
As they best exploit the structure of the suboperator, we would expect FCG and PIPEFCG to be more performant.
Table \ref{tab:ex43} shows the resulting timings. 
We see significant speedup in all cases as the reduction bottlenecks in the base solvers are overcome.
In the case where Conjugate Residuals are employed as inner solvers, we observe variability in the 
number of iterations to convergence, due to the degradation in conditioning arising from use of the 
$A^HA$-inner product.

\begin{table}[htbp]
\begin{center}
\footnotesize
\begin{tabular}{lll ll}
\toprule
               &   & \multicolumn{3}{c}{Solve Time  [s]} \\
  Inner Solver &   Outer Its. & Mean (Std. Dev.)        & Min.     & Max.  \\
\midrule
  FCG          & 8            & 3.80 (0.29)         & 3.47    & 4.37  \\
  PIPEFCG  & 8 & 1.80 (0.18) & 1.52  & 2.03     \\ 
\midrule
  GCR & 18--23 & 10.02 (1.51) & 8.26 & 12.76    \\
  PIPEGCR & 12--16 & 3.08 (0.46) & 2.59 & 4.01  \\
\midrule
  FGMRES & 19 & 6.70 (0.53) & 6.10 & 7.71      \\
  PIPEFGMRES & 19 & 5.04 (0.27) & 4.73 & 5.65  \\
\bottomrule
\end{tabular}
%\begin{tabular}{lcccccccccc}
%\hline
%               &   & \multicolumn{4}{c}{Linear Solve Time / Iteration [s]}\\
%  Inner Solver &   \white{Outer Its.}  &  Mean (std. dev.)        & Min      & Median   & Max  \\
%\hline
%  FCG        & \white{8    }& 3.80 (0.29) & 3.47 & 3.70 & 4.37 \\
%  PIPEFCG    & \white{8    }& 0.23 (0.02) & 0.19 & 0.23 & 0.25 \\ 
%  GCR        & \white{18-23}& 0.47 (0.04) & 0.42 & 0.46 & 0.55 \\
%  PIPEGCR    & \white{12-16}& 0.24 (0.02) & 0.22 & 0.25 & 0.25 \\
%  FGMRES     & \white{19   }& 0.35 (0.03) & 0.32 & 0.34 & 0.41 \\
%  PIPEFGMRES & \white{19   }& 0.27 (0.01) & 0.25 & 0.26 & 0.30 \\
%\hline
%\end{tabular}
\caption{Statistics for 9 runs of a 2D Stokes flow problem with an upper block triangular preconditioner employing various inner solvers, as described in \S\ref{sec:ex43} \label{tab:ex43}.}
\end{center}
\end{table}

%%%%%%%%%%%%%%%%%%%%%%%%%%%%%%%%%%%%%%%%%%%%%%%%%%%%%%%%%%%%%%%%%%%%%%%%%%%%%
%%%%%%%%%%%%%%%%%%%%%%%%%%%%%%%%%%%%%%%%%%%%%%%%%%%%%%%%%%%%%%%%%%%%%%%%%%%%%
\section{Performance Model and Extrapolation to Exascale}
%%%%%%%%%%%%%%%%%%%%%%%%%%%%%%%%%%%%%%%%%%%%%%%%%%%%%%%%%%%%%%%%%%%%%%%%%%%%%
%%%%%%%%%%%%%%%%%%%%%%%%%%%%%%%%%%%%%%%%%%%%%%%%%%%%%%%%%%%%%%%%%%%%%%%%%%%%%
\label{sec:perf-model}
\label{sec:model-exasc-mach}

We describe performance models to predict the performance of our methods
at exascale. 
For ease of comparison we consider the same hypothetical exaflop
machine used to analyze pipelined GMRES \cite{Ghysels2013}, based on a 2010
report \cite{Shalf2010}. 
The specifications of the hypothetical machine are
listed in Table~\ref{tab:examachine_specs}. 
Rank-local computation times are
obtained by counting flops and multiplying by \tc.
\revOne{7}{This makes the significant assumption that any give computation can attain the peak flop rate of the machine, without being memory bandwidth limited.}

\begin{table}[h!]
\footnotesize %from SIAM latex guidelines
  \begin{center}
    \begin{tabular}{lll}
      \toprule
      Property & Symbol & Value\\
      \midrule
      Nodes & \Pn & $2^{20}$\\
      Cores per node & $C_\text{n}$ & $2^{10}$\\
      Compute cores& $\Pc = \Pn C_\text{n}$ & $2^{30}$\\
      Word size &   $w$      & 32 B\\
      Machine node interconnect bandwidth & BW  & 100 GB/s\\
      Tree radix & r            & 8\\
      Per word transfer time & \tw           & $w$/BW\\
      Latency (startup) time & \ts           & 1 $\mu$s        \\
      Time per flop & \tc           & $2^{30}/10^{18}$  \\
      \bottomrule
    \end{tabular}
    \caption{Properties of a hypothesized exascale machine.}
    \label{tab:examachine_specs}
  \end{center}
\end{table}
 
\paragraph{Global Reductions}
Dot products are typically computed with a reduction tree of height 
$\lceil \log_r(\Pn)\rceil$ with $r$ the tree radix.
Present day networks make use of a hierarchy of interconnects where
more closely located nodes, e.g.~in a cabinet, are connected through a tree with
a high effective radix. 
A job scheduler typically assigns compute nodes, giving the user 
limited influence on the specific compute nodes used.
Our performance model makes the simplifying assumption of a tree radix $r = 8$ for
the entire machine \cite{Ghysels2013,Shalf2010}.

In order to obtain an expression for the global reduction communication
  time $\Tredcomm$ we assume the model $\ts + m\tw$ consisting of a constant
latency or startup time \ts\ and a message size dependent part involving the
number of words $m$ and the time \tw\ required to exchange a word. 
We neglect time for intranode communication leading to the cost of $2\lceil \log_r(\Pn)\rceil
(\ts + m\tw)$ for a reduction and subsequent broadcast across the entire machine.
In the context of pipelined methods global reduction costs include only the
excess time (if any) which is not overlapped by other work.

Global reductions also involve floating point operations, accounted for in the 
computation time \Tredcalc. Adding up the partial reduction results
across the reduction tree leads to $\lceil \log_r (P_n )\rceil$ flops along the
longest path. While this is the maximum number of flops performed by one
process only, all other processes wait for this summation to finish. 
This might appear to be an insignificant contribution to the overall
computational cost but can become non-negligible in the strong scaling limit.

\paragraph{Sparse Matrix-Vector Products (SpMV)}

The cost of a sparse matrix-vector multiplication also includes both
communication and computation. 
The system matrix is assumed to represent
a stencil based approximation of the underlying partial differential
equations, with $\nz$ \revOne{25}{non-zero} entries per row.
For a standard second order finite difference discretization in three dimensional
space, e.g., $\nz = 7$. 
Denoting by $N = N_x N_y N_z$ the total problem size and by $\nloc = N/\Pc$ the portion local
to a process we obtain the compute costs $\Tspmvcalc = 2\nz\nloc\tc$. 
Neighbor communication costs are $\Tspmvcomm =  6(t_s +  (N/\Pn)^{\nicefrac{2}{3}}t_w)$ and can
mostly be overlapped. 
Consequently, the SpMV cost function is $\Tspmv = \max(\Tspmvcalc,\Tspmvcomm)$.

\begin{table}[h!tb]
\footnotesize %from SIAM latex guidelines
  \begin{center}
    \begin{tabular}{lll}
      \toprule
      Method  & Operation & Communication time cost function \\
      \midrule
      RED & for $k = 1..w:$ $\beta_k \gets \langle p,u_k \rangle$ & {$\!\begin{aligned}
          \Tredcomm(w) &= 2 \lceil \log_r (P_n )\rceil \big(t_s + w t_w\big)
        \end{aligned}$}\\
      SpMV & $s \gets A p $ & {$\!\begin{aligned}
          \Tspmvcomm &= 6(t_s +  (N/P_n)^{\nicefrac{2}{3}}t_w)
        \end{aligned}$}\\
      \bottomrule
    \end{tabular}
    \caption{Communication time cost functions for fundamental operations with $\nloc$
      as defined above, $\nz$ the number of nonzero entries per row, and $w$ the
      number of values (words) to be broadcast. }
    \label{tab:performance_fundamentals_comm}
    \revOne{26}{}
  \end{center}
\end{table}

\paragraph{Preconditioner Application}

We augment the model with the application of
a preconditioner.
\revOne{6}{
We model the performance of the RASM preconditioner used in} \S\ref{sec:tests}, \revOneCol{defining $\Tpc$ to account for 5 CG-Jacobi iterations on subdomains overlapping by one degree of freedom, and an internode communication step with the same pattern as a SpMV.}
\revOneCol{We note that scalable multigrid preconditioners may also be designed which do not involve any reductions} \cite{May:2015iu} \revOneCol{and thus could be attractive for extreme scale solves. Nested pipelined Krylov methods could be used as coarse grid solvers within these preconditioners.}

\paragraph{Total Cost per Iteration}

The total cost of one iteration comprises the time for performing
flops, $\Tcalc$, involving the cost of sparse matrix-vector
products, $\Tspmv$, the time for applying the preconditioner, $\Tpc$, and fundamental
operations listed in Table~\ref{tab:performance_fundamentals} as well as the (excess of
the) reduction communication time $\Tredcomm$.
We do not take into account fused operations and vector units
to establish an easier comparison with previous work \cite{Ghysels2013}.

The cost per iteration for flexible methods depends on the
number of previous directions to orthogonalize against and hence on the
truncation strategy and the iteration count.
In our model we use a factor
representing an average number of previous directions being used in the
Gram-Schmidt orthogonalization as $\nuavg=\kavg \cdot \numax$.
Operations involving scalar values only are not taken into account. This leads
to the models for the single iteration cost functions listed in Table~\ref{tab:performance_model}.

\begin{table}[htbp]
\footnotesize %from SIAM latex guidelines
  \begin{center}
    \begin{tabular}{lll}
      \toprule
      Method & Operation & Compute time cost function \\
      \midrule
      AXPY & $x \gets x + \alpha p$ & {$\!\begin{aligned}
          \Taxpy &= t(\nloc \text{(\texttt{mlt})} + \nloc \text{(\texttt{add})}) = 2 \nloc t_c
        \end{aligned}$}\\
      MAXPY & $p \gets u + \sum_m \beta_m p_m$ & {$\!\begin{aligned}
          \Tmaxpy(m) &= t(m \nloc \text{(\texttt{mlt})} + (m-1) \nloc \text{(\texttt{add})} \\
          & +\nloc \text{(\texttt{add})}) = 2 m \nloc t_c
        \end{aligned}$}\\
      RED & $\langle r,u \rangle$ & {$\!\begin{aligned}
          \Tredcalc &= t(\nloc \text{(\texttt{mlt})} + (\nloc-1) \text{(\texttt{add})} \\ 
          &+  \lceil \log_r (P_n)\rceil \text{(\texttt{add})}) \approx (2 \nloc + \lceil \log_r (P_n )\rceil) t_c
        \end{aligned}$}\\
      \bottomrule
    \end{tabular}
    \caption{Compute time cost functions for fundamental operations with $\nloc=N/P_c$
      denoting the vector size local to the process.}
    \label{tab:performance_fundamentals}
  \end{center}
\end{table}

\begin{table}[htbp]
\footnotesize %from SIAM latex guidelines
  \begin{center}
    \begin{tabular}{l|l}
      \toprule
      Krylov method &Total cost \\
	\midrule      
      FCG & {$\!\begin{aligned}
          \Tcalc &= 2 \Taxpy + \Tmaxpy(\nuavg) + (\nuavg+2)\Tredcalc + \Tspmv + \Tpc\\
          \Tred &= \Tredcomm(\nuavg+1) + \Tredcomm(1)
        \end{aligned}$}
      \\
	\midrule      
      PIPEFCG & {$\!\begin{aligned}
          \Tcalc &= 4 \Taxpy + 4 \Tmaxpy(\nuavg) + (\nuavg+2)\Tredcalc + \Tspmv \\ 
          &+ \Tpc + 2 n t_c\\
          \Tred &= \max\big(0,\Tredcomm(\nuavg+2) - [\Tspmv+\Tpc+2 n t_c] \big)
        \end{aligned}$}
      \\
	\midrule      
      GCR & {$\!\begin{aligned} 
          \Tcalc &= 2 \Taxpy + \Tmaxpy(\nuavg) + (\nuavg+2)\Tredcalc + \Tspmv + \Tpc\\
          \Tred &= \Tredcomm(\nuavg) + \Tredcomm(2)
        \end{aligned}$} \\
	\midrule      
      PIPEGCR & {$\!\begin{aligned} 
          \Tcalc &= 3 \Taxpy + 3 \Tmaxpy(\nuavg) + (\nuavg+2)\Tredcalc + \Tspmv \\ 
          &+ \Tpc+2 n t_c\\
          \Tred &= \max\big(0,\Tredcomm(\nuavg+2) - [\Tpc + 2 n t_c] \big)
        \end{aligned}$} \\
	\midrule      
      PIPEGCR\_w & {$\!\begin{aligned} 
          \Tcalc &= 4 \Taxpy + 4 \Tmaxpy(\nuavg) + (\nuavg+2)\Tredcalc + \Tspmv \\
          &+ \Tpc+2 n t_c\\
          \Tred &= \max\big(0,\Tredcomm(\nuavg+2) - [\Tspmv+\Tpc+2 n t_c] \big)
        \end{aligned}$} \\
	\midrule      
      FGMRES & {$\!\begin{aligned} 
          \Tcalc &= \Tmaxpy(\nuavg) + t_c n + (\nuavg+1)\Tredcalc + \Tspmv + \Tpc\\
          \Tred &= \Tredcomm(\nuavg) + \Tredcomm(1)
        \end{aligned}$} \\
	\midrule      
      PIPEFGMRES & {$\!\begin{aligned} 
          \Tcalc &= \nuavg \Taxpy + 3 \Tmaxpy(\nuavg) + (\nuavg+5) t_c n  + (\nuavg+2)\Tredcalc\\                  &\quad + \Tspmv + \Tpc\\
          \Tred &= \max\big(0,\Tredcomm(\nuavg+2) - [\Tspmv + \Tpc] \big)
        \end{aligned}$} \\
      \bottomrule
    \end{tabular}
    \caption{Performance models for standard and pipelined methods for the
      truncation strategy. Again, $n = N/P_c$, $\nuavg=\kavg\cdot\numax$, where $\kavg$ corresponds to
      a factor accounting for \revOneCol{the average number of previous search directions being used}. 
      If the number of iterations required is $\gg \numax$, this factor
      approaches one for a standard truncation strategy. 
      Communication time for multiple reductions is fused whenever
      possible. GMRES-type models do not include the time for
      manipulating the Hessenberg matrix. 
      }
    \label{tab:performance_model}
  \end{center}
\revOne{10}{} %this is a hack, but we don't expect many of these comments inside figures
\end{table}

\paragraph{Problem and Solver Specifications}

For evaluating the cost functions, we consider a strong scaling experiment (cf. \cite{Ghysels2013}) with
a total problem size $N = 2000^3$, $\nz = 7$ non-zero entries per
row and evaluate for a varying number of node counts. For {PIPEFCG} and
{PIPEGCR} we assume $\numax = 30$ and
$\kavg = 0.8$; for {PIPEFGMRES} we assume a restart parameter of 30. 

\paragraph{Evaluation}

\begin{figure}
  \revOne{6}{}
  \begin{center}
    \includegraphics[trim=0.75cm 0.5cm 1.75cm 1cm,width=0.48\textwidth]{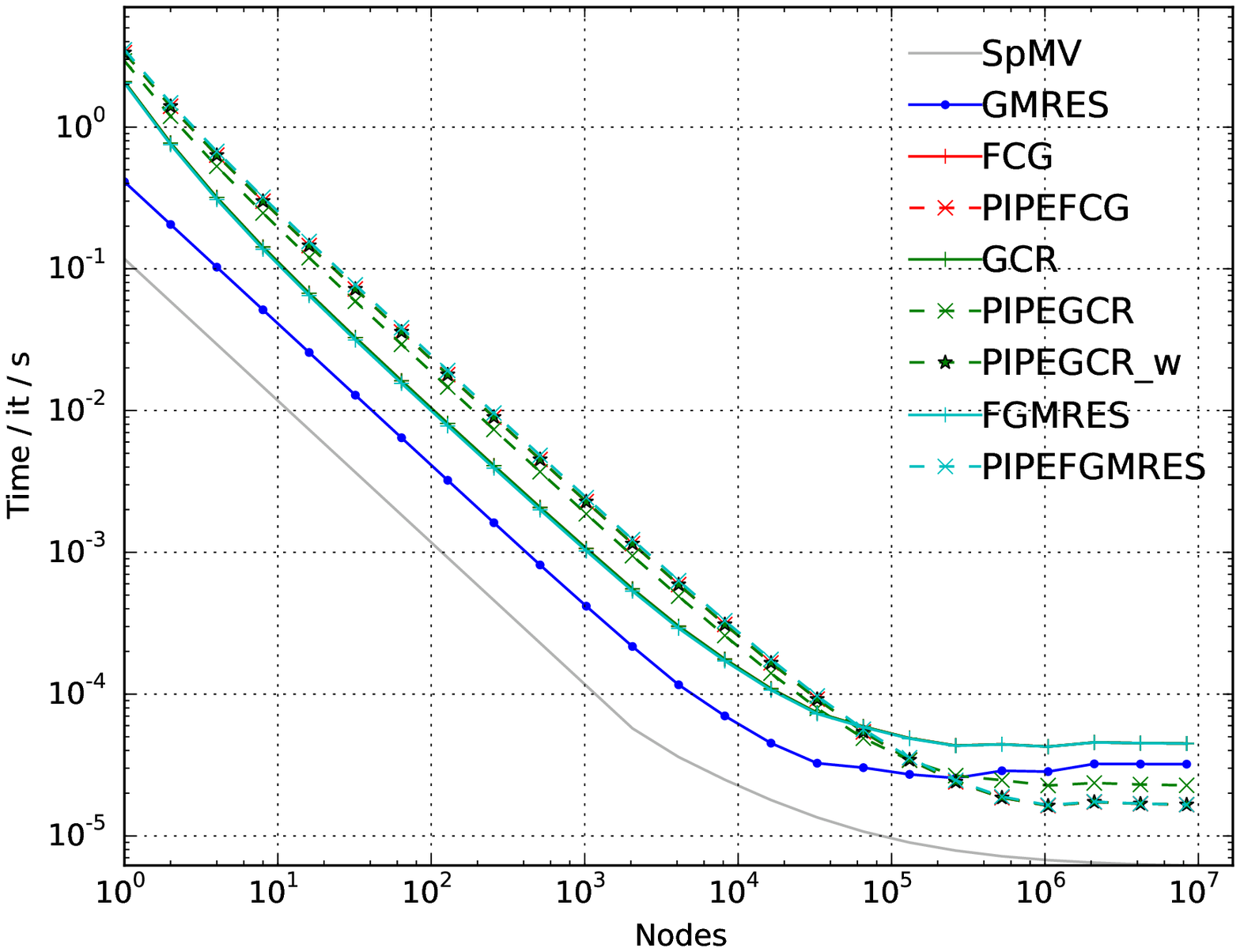}\hspace{2mm}
    \includegraphics[trim=0.75cm 0.5cm 1.75cm 1cm,width=0.48\textwidth]{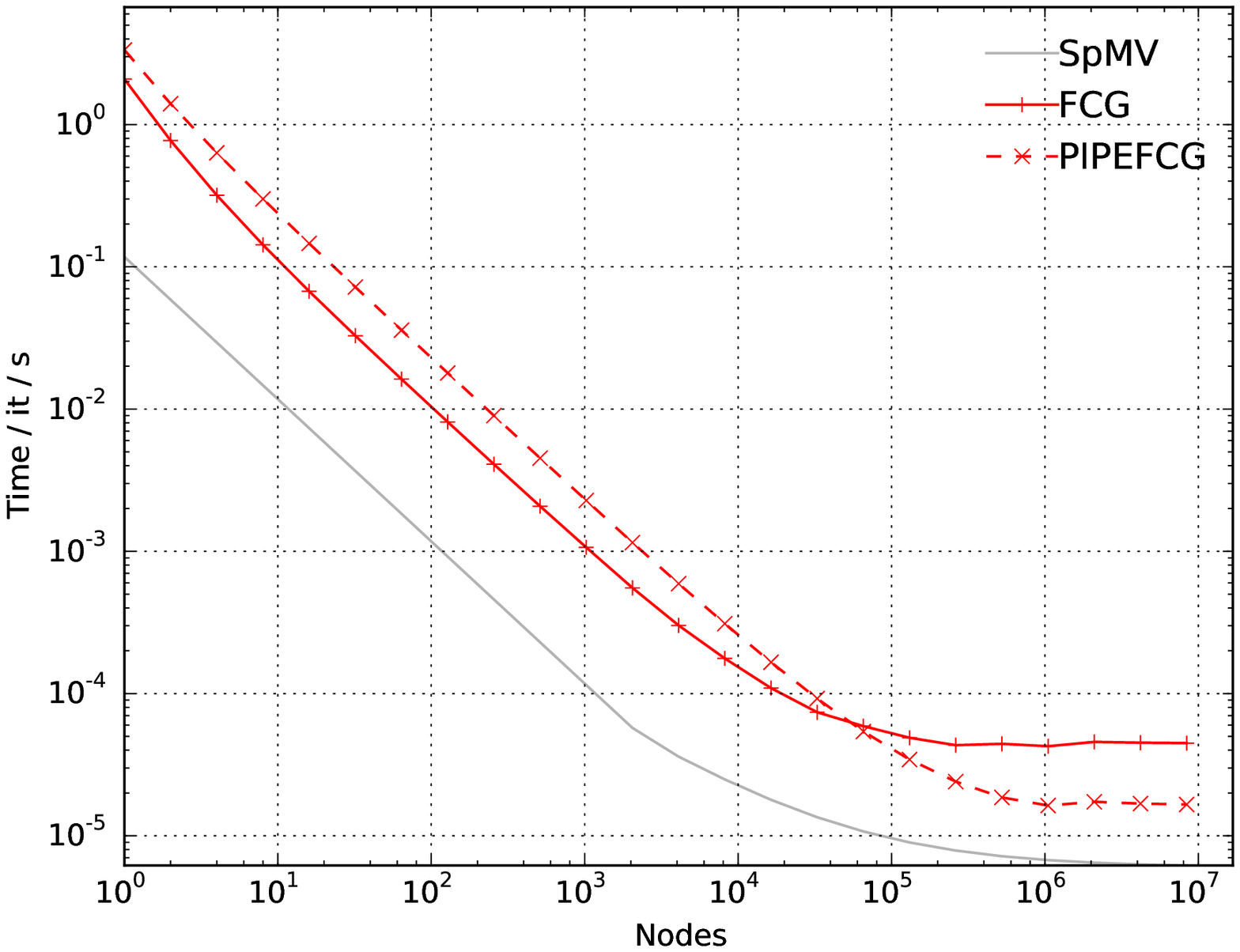}\\~\\
    \includegraphics[trim=0.75cm 0.5cm 1.75cm 1cm,width=0.48\textwidth]{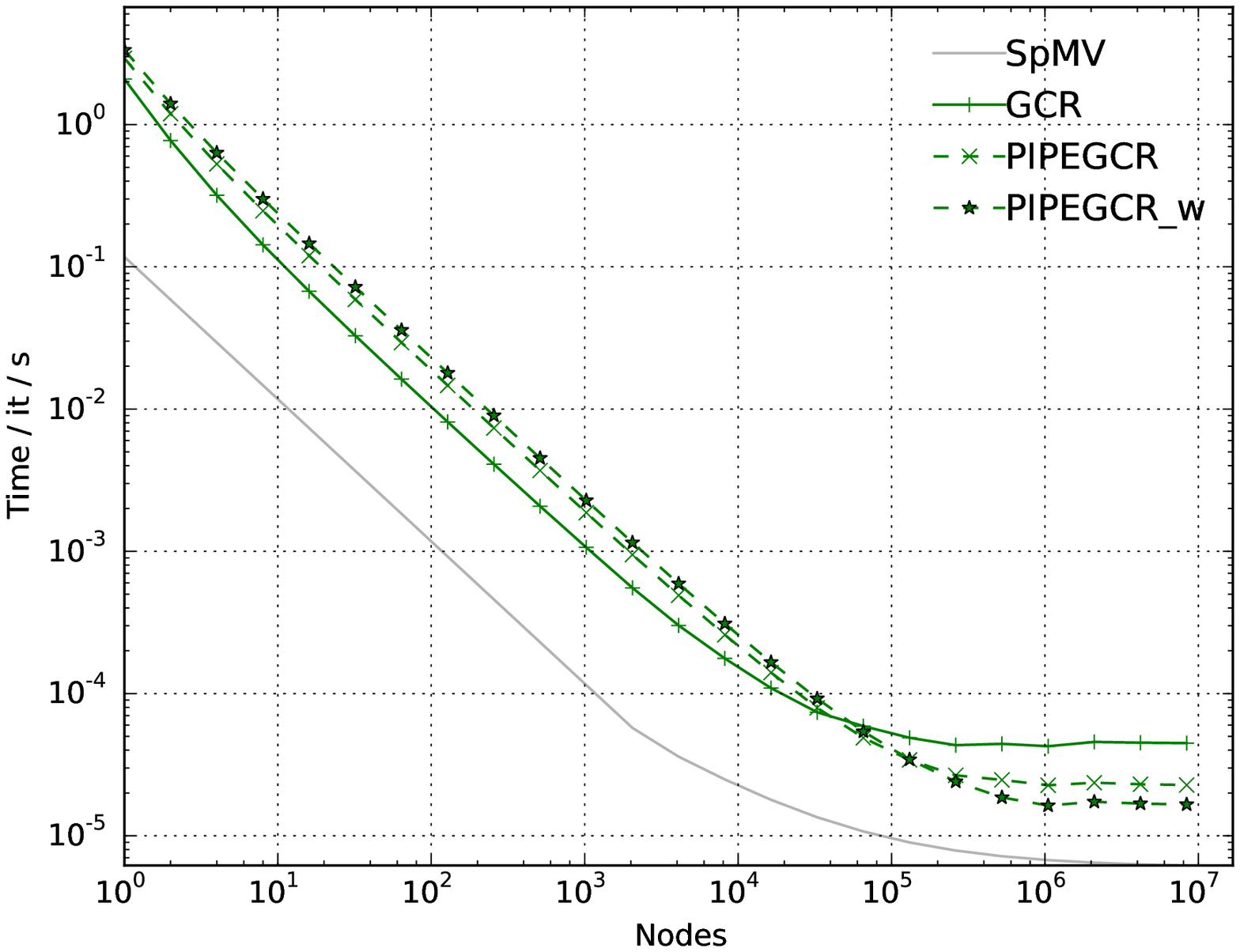}\hspace{2mm}
    \includegraphics[trim=0.75cm 0.5cm 1.75cm 1cm,width=0.48\textwidth]{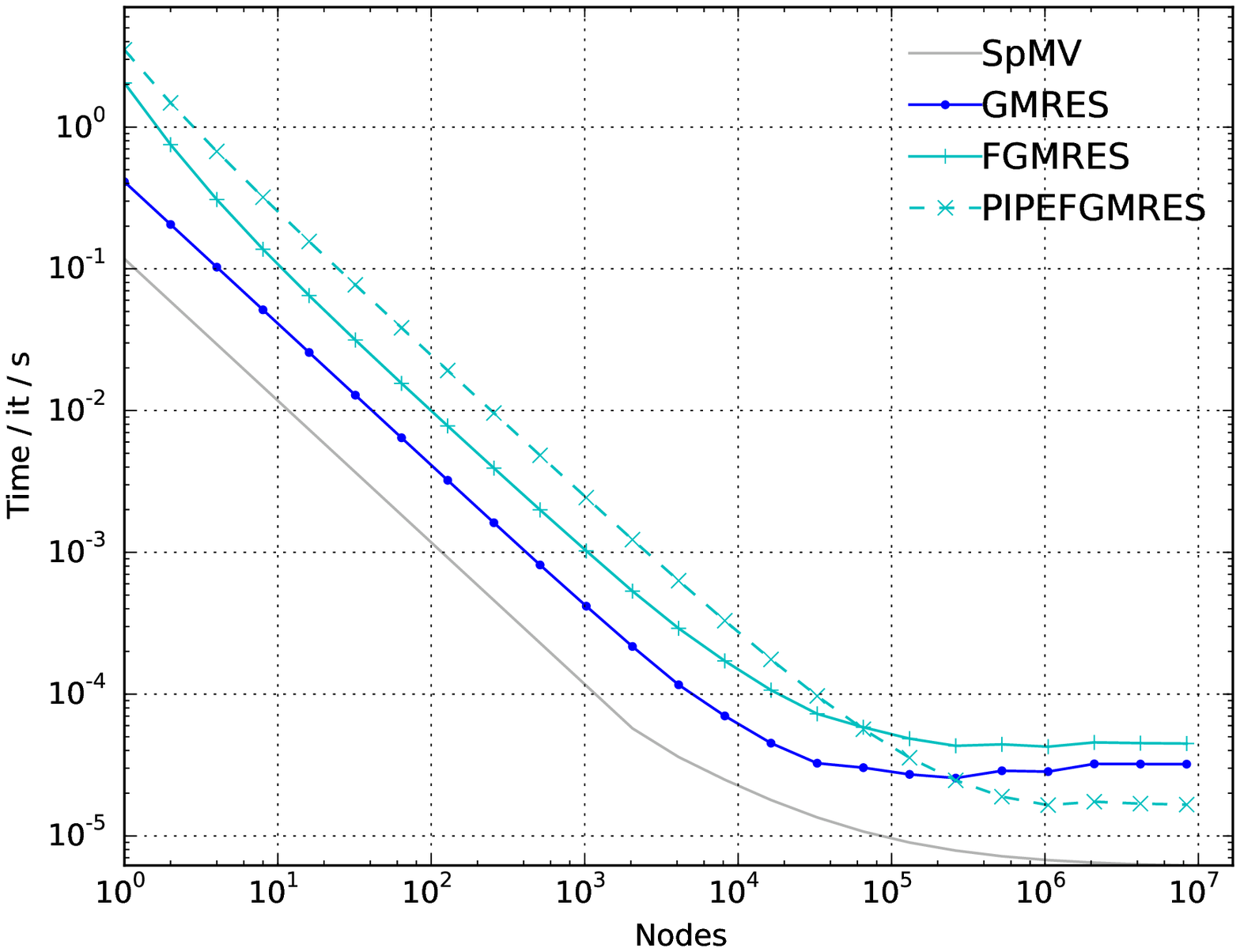}
  \end{center}
  \caption{Evaluation of the performance models of
    Table~\ref{tab:performance_model} for predicting iteration times on the
    hypothesized exascale machine of Table~\ref{tab:examachine_specs}. The top
    left panel includes all methods and the others compare FCG-, GCR- and
    GMRES-type methods. All graphs include the time for performing a {SpMV} operation
    as a reference. The number of nodes scales beyond the 1M nodes of
    the assumed machine to better visualize the leveling at high node
    counts. While all methods are eventually limited by latency time the
    pipelined methods outperform their standard counterparts on high node counts.}
  \label{fig:performance_model_comp}
\end{figure}

In Figure~\ref{fig:performance_model_comp} we show a comparison of the
anticipated performance of all solvers for
this problem. 
Eventually all methods are limited by latency time but the
graphs indicate that the pipelined methods scale up to the $10^6$ nodes
of the assumed exascale machine while
their standard counterparts level off about one order of magnitude earlier. The
performance crossover regarding pipelined vs.~non-pipelined methods for this
particular example occurs around \revOne{6}{$8\cdot 10^4$} nodes. The example of
{PIPEGCR} vs.~{PIPEGCR\_w} clearly shows the effect
of trading local work load for better overlapping. Note how {PIPEGCR}
(explicit computation of $w$, overlapping by preconditioner application)
performs slightly better at lower node counts while {PIPEGCR\_w} (recurrence of $w$, additional
{AXPY} and {MAXPY}, overlapping by preconditioner and
{SpMV}) shows a more significant gain in efficiency for larger numbers of nodes.

\begin{figure}
  \begin{center}
    \includegraphics[trim=0.75cm 0.5cm 1.75cm 1cm,width=0.48\textwidth]{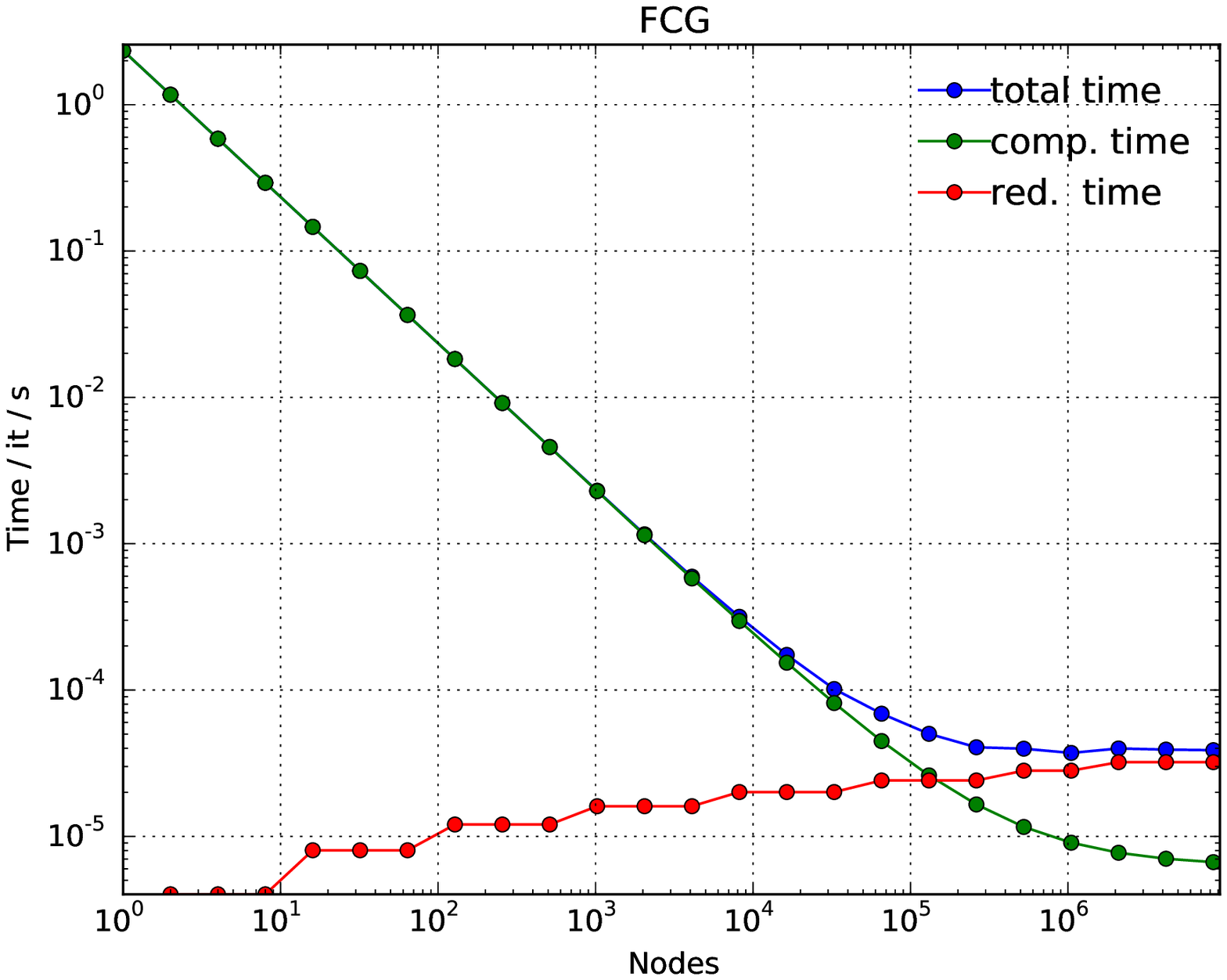}\hspace{2mm}
    \includegraphics[trim=0.75cm 0.5cm 1.75cm 1cm,width=0.48\textwidth]{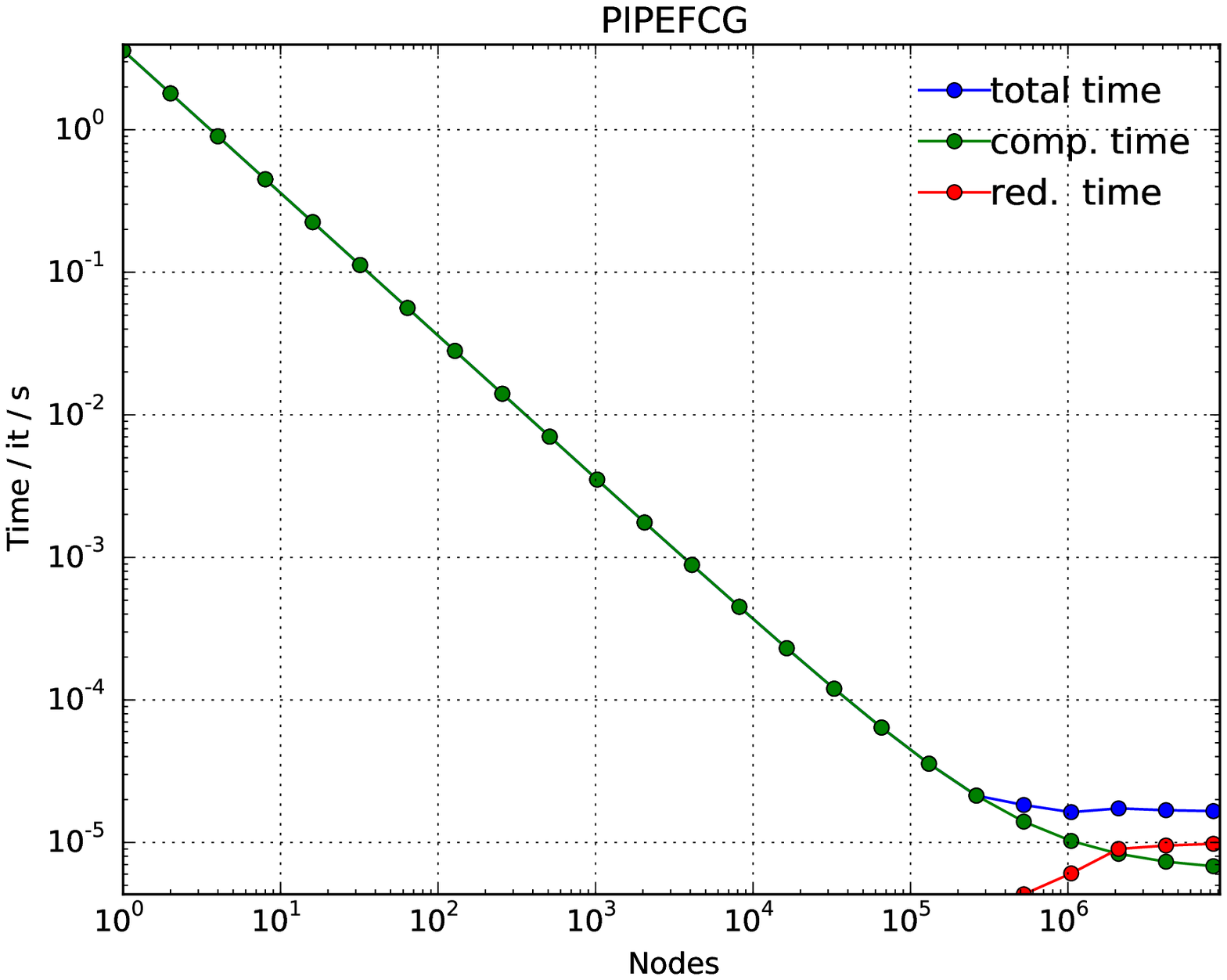}
  \end{center}
  \caption{Breakdown of the total iteration time into computation and reduction
    time for FCG and PIPEFCG. For the test case considered it is anticipated
    that the reductions will be
    entirely overlapped up to approximately $5\cdot 10^5$ nodes for the pipelined method.}
  \label{fig:performance_model_breakdown}
\end{figure}

Figure~\ref{fig:performance_model_breakdown} shows a breakdown of the
anticipated total iteration time for {FCG} and {PIPEFCG}. For the
test case considered here the model predicts that global reductions will be
entirely overlapped up to approximately $5\cdot 10^5$ nodes.

\paragraph{Sensitivity to Parameters}

The performance models are insensitive to the specific
value of most of the parameters when chosen within reasonable bounds. Among the insensitive
parameters are the interconnect bandwidth, word size, the tree radix and
assumptions such as ignoring fused {AXPY} operations.
It may appear unintuitive that the interconnect bandwidth is not a sensitive
parameter, but latency is dominant in the communication time model.

The models are most obviously sensitive to the latency time but also to
the number of cores per node as more or less cores per node implies a different
number of nodes for obtaining exaflop performance. Another sensitive parameter
is the total problem size $N$. A mild dependency of the average fill factor \kavg\ and the
maximum number of previous direction \numax\ as well as the number of non-zero
entries \nz\ per row is observed. Varying any of these values does not cause a change
of the general shape of the curves but rather induces a shift of the crossover point.

\revOne{28}{For the} hypothesized exascale machine a latency
time of $1 \mu s$ is assumed, significantly lower than times observed on current
systems including Piz Daint as in Table~\ref{tab:reduction_timing}. 
Consequently, pipelined
Krylov methods can lead to significant performance improvements for much smaller
numbers of nodes on current systems, as confirmed by the tests in \S\ref{sec:tests}.

%%%%%%%%%%%%%%%%%%%%%%%%%%%%%%%%%%%%%%%%%%%%%%%%%%%%%%%%%%%%%%%%%%%%%%%%%%%%%
%%%%%%%%%%%%%%%%%%%%%%%%%%%%%%%%%%%%%%%%%%%%%%%%%%%%%%%%%%%%%%%%%%%%%%%%%%%%%
\section{Conclusions and Outlook}
%%%%%%%%%%%%%%%%%%%%%%%%%%%%%%%%%%%%%%%%%%%%%%%%%%%%%%%%%%%%%%%%%%%%%%%%%%%%%
%%%%%%%%%%%%%%%%%%%%%%%%%%%%%%%%%%%%%%%%%%%%%%%%%%%%%%%%%%%%%%%%%%%%%%%%%%%%%

Bottlenecks related to global communication and synchronization are increasingly
common in high performance computing, as parallelism and hybridism are  used to 
increase peak performance.
Traditional algorithms can often be improved with this consideration in mind.
Given the increasing relative impact of latency and synchronization required with global all-to-all communication patterns, 
it becomes increasingly attractive to tailor algorithms to mitigate these costs.
This requires introducing additional costs in terms of memory footprint, local computational work, memory traffic, numerical stability, and allowing asynchronous operations.

We introduce variants of the FCG, GCR and FGMRES flexible Krylov methods, allowing 
overlap of global reductions with other work
through operation pipelining. 
The base methods are amongst the most commonly used 
flexible Krylov methods for Hermitian positive definite, Hermitian, and general
operators, respectively.
Future work includes extending and implementing the approaches here
for other useful Krylov methods.
A crucial consideration in this process, as demonstrated here, is the effect of 
allowing algorithmic rearrangements which modify the effective preconditioner.
Naive rearrangement is not always effective, and effective rearrangements may
impose conditions on the class of useful preconditioners beyond those common
in previous Krylov methods.
We propose and analyze an approach to provide methods which are increasingly
effective for stronger preconditioners.

The proposed methods are shown to be effective on a current leadership supercomputing platform with a fast network, where speedups could exceed $2\times$, challenging common assumptions that pipelined Krylov methods are only effective at future machine scales, and that the speedup of single-stage pipelining is limited to $2\times$.
Additionally, we develop analytic
models for anticipating performance on future exascale
systems.

All methods have been implemented in the \textsc{PETSc} package and are available open-source.

\paragraph{Deeper Pipelining} The same rearrangements used in the single-stage \\
pipelined methods described here could be pursued to allow a reduction to overlap more 
 than one iteration's worth of other work. 
 This has been explored for the GMRES method \cite{Ghysels2013}. 
Our focus has been on the ``faithful'' preconditioners discussed above, which require
 a substantial amount of work to apply and thus would likely not benefit from further pipelining.
Additionally, numerical instability becomes more of a concern with deeper pipelining 
--in the case of GMRES, \revOne{2}{norm} breakdown becomes more common--and we anticipate 
that the required stabilization for nonlinear preconditioners would also be non-trivial. 
Investigation of deeper pipelining is of course of interest from a mathematical 
point of view, but it is unclear whether exascale systems will be developed
exhibiting the extreme reduction latency required to justify overlapping inner 
products with multiple nonlinear preconditioner applications. 
From a practical point of view, current MPI implementations do not prioritize 
operations like multiple overlapping reductions, which makes the testing of these methods problematic. 

\paragraph{Outlook}
The methods described in this paper are expected to become more and more relevant as non-deterministic, randomized, nested, highly-distributed, and nonlinear preconditioning techniques come into greater usage. 
These trends are encouraged by an increasing level of parallelism, hybridism, and hierarchy in computational machinery, as well as intense research into randomized, finely-grained parallel, and asynchronous techniques for approximate and exact solvers.

\section*{Acknowledgments}
We acknowledge financial support from the Swiss University Conference and the Swiss Council of Federal Institutes of Technology through the Platform for Advanced Scientific Computing (PASC) program.

\ifx\undefined\nosiam
  \bibliographystyle{siamplain}
\else
  \bibliographystyle{siam}
\fi
\bibliography{PFKrylov}

\end{document}